\documentclass[11pt]{amsart}
\usepackage{latexsym,graphicx}
\numberwithin{equation}{section}
\theoremstyle{plain}

\theoremstyle{remark}

\theoremstyle{definition}

\newcommand{\D}{{\mathcal D}}
\newcommand{\E}{\mathcal E}

\newcommand{\G}{{\mathcal G}}

\newcommand{\K}{{\mathcal K}}
\renewcommand{\L}{{\mathcal L}}
\newcommand{\M}{{\mathcal M}}
\newcommand{\N}{\mathbb N}

\newcommand{\V}{{\mathcal V}}

\newcommand{\dist}{\operatorname{dist}}

\newcommand{\esssup}{\operatorname{ess\,sup}}
\newcommand{\fp}{\operatorname{FP}}

\newcommand{\Int}{\operatorname{Int}}

\renewcommand{\span}{\operatorname{span}}

\newcommand{\supp}{\operatorname{Supp}}

\def\half{{1 \over 2}}

\newcommand{\oa}{\overrightarrow}

\newcommand{\ol}{\overline}

\def\XXint#1#2#3{{\setbox0=\hbox{$#1{#2#3}{\int}$}
      \vcenter{\hbox{$#2#3$}}\kern-.5\wd0}}

\begin{document}

\def\cal{\mathcal}

\font\tpt=cmr10 at 12 pt
\font\fpt=cmr10 at 14 pt

\font \fr = eufm10

%\font\AAA=Times.dfont  at 12pt
 %\font\BBB=Times.dfont at 8pt

%\font\AAA=cmr10 at 12pt
%\font\BBB=cmr10 at 8pt

\def\AAA{\bf}
\def\BBB{\bf}

\overfullrule=0in

\def\boxit#1{\hbox{\vrule
 \vtop{%
  \vbox{\hrule\kern 2pt %
     \hbox{\kern 2pt #1\kern 2pt}}%
   \kern 2pt \hrule }%
  \vrule}}

  \def\harr#1#2{\ \smash{\mathop{\hbox to .3in{\rightarrowfill}}\limits^{\scriptstyle#1}_{\scriptstyle#2}}\ }

\def\ALEX{1}
\def\AGV{2}
\def\ASSA{3}
\def\ASSB{4}
\def\BTA{5}
\def\BTB{6}
\def\BTC{7}
\def\CLN{8}
\def\CIL{9}
\def\CRA{10}
\def\GAR{11}
\def\HAR{12}
\def\DDD{13}
\def\PUP{14}
\def\DDR{15}
\def\HYP{16}
\def\HLGAR{17}
\def\REST{18}
\def\AC{19}
\def\pPSH{20}
\def\SURVEY{21}
\def\BELL{22}
\def\ASPECTS{23}
\def\HP{24}
\def\HPP{25}
\def\LAB{26}
\def\LLAB{27}
\def\LLLAB{28}
\def\LAN{29}
\def\LEL{30}
\def\POGA{31}
\def\POGB{32}
\def\SHF{33}
\def\TWA{34}
\def\TWB{35}
\def\TWC{36}
\def\VIT{37}

 \def\GG{{{\bf G} \!\!\!\! {\rm l}}\ }

\def\GL{{\rm GL}}

\def\bll{I \!\! L}

\def\bra#1#2{\langle #1, #2\rangle}
\def\bbf{{\bf F}}
\def\bbj{{\bf J}}
\def\Jtn{{\bbj}^2_n}  \def\JtN{{\bbj}^2_N}  \def\JoN{{\bbj}^1_N}
\def\jt{j^2}
\def\jtx{\jt_x}
\def\Jt{J^2}
\def\Jtx{\Jt_x}
\def\bpp{{\bf P}^+}
\def\bpt{{\wt{\bf P}}}
\def\fsh{$F$-subharmonic }
\def\mo{monotonicity }
\def\jet{(r,p,A)}
\def\ss{\subset}
\def\sse{\subseteq}
\def\half{\hbox{${1\over 2}$}}
\def\smfrac#1#2{\hbox{${#1\over #2}$}}
\def\oa#1{\overrightarrow #1}
\def\dim{{\rm dim}}
\def\dist{{\rm dist}}
\def\codim{{\rm codim}}
\def\deg{{\rm deg}}
\def\rank{{\rm rank}}
\def\log{{\rm log}}
\def\Hess{{\rm Hess}}
\def\Hessyp{{\rm Hess}_{\rm SYP}}
\def\trace{{\rm trace}}
\def\tr{{\rm tr}}
\def\max{{\rm max}}
\def\min{{\rm min}}
\def\span{{\rm span\,}}
\def\Hom{{\rm Hom\,}}
\def\det{{\rm det}}
\def\End{{\rm End}}
\def\Sym{{\rm Sym}^2}
\def\diag{{\rm diag}}
\def\pt{{\rm pt}}
\def\Spec{{\rm Spec}}
\def\pr{{\rm pr}}
\def\Id{{\rm Id}}
\def\Grass{{\rm Grass}}
\def\Herm#1{{\rm Herm}_{#1}(V)}
\def\arr{\longrightarrow}
\def\supp{{\rm supp}}
\def\Link{{\rm Link}}
\def\Wind{{\rm Wind}}
\def\Div{{\rm Div}}
\def\vol{{\rm vol}}
\def\foral{\qquad {\rm for\ all\ \ }}
\def\fpsh{{\cal PSH}(X,\f)}
\def\Core{{\rm Core}}
\def\dis{f_M}
\def\Re{{\rm Re}}
\def\rn{\bbr^n}
\def\pp{\cp^+}
\def\plp{\cp_+}
\def\Int{{\rm Int}}
\def\cix{C^{\infty}(X)}
\def\Gr#1{G(#1,\rn)}
\def\Symn{{\Sym(\rn)}}
\def\SymN{{\Sym(\bbr^N)}}
\def\Gpn{G(p,\rn)}
\def\fd{{\rm free-dim}}
\def\SA{{\rm SA}}
 \def\cd{{\cal C}}
 \def\cdt{{\widetilde \cd}}
 \def\cm{{\cal M}}
 \def\cmt{{\widetilde \cm}}

\def\Theorem#1{\medskip\noindent {\bf THEOREM \bf #1.}}
\def\Prop#1{\medskip\noindent {\bf Proposition #1.}}
\def\Cor#1{\medskip\noindent {\bf Corollary #1.}}
\def\Lemma#1{\medskip\noindent {\bf Lemma #1.}}
\def\Remark#1{\medskip\noindent {\bf Remark #1.}}
\def\Note#1{\medskip\noindent {\bf Note #1.}}
\def\Def#1{\medskip\noindent {\bf Definition #1.}}
\def\Claim#1{\medskip\noindent {\bf Claim #1.}}
\def\Conj#1{\medskip\noindent {\bf Conjecture \bf    #1.}}
\def\Ex#1{\medskip\noindent {\bf Example \bf    #1.}}
\def\Qu#1{\medskip\noindent {\bf Question \bf    #1.}}
\def\Exercise#1{\medskip\noindent {\bf Exercise \bf    #1.}}

\def\HoQu#1{ {\AAA T\BBB HE\ \AAA H\BBB ODGE\ \AAA Q\BBB UESTION \bf    #1.}}

\def\pf{\medskip\noindent {\bf Proof.}\ }
\def\qed{\hfill  $\vrule width5pt height5pt depth0pt$}
\def\equdef{\buildrel {\rm def} \over  =}
\def\qedqed{\hfill  $\vrule width5pt height5pt depth0pt$ $\vrule width5pt height5pt depth0pt$}
\def\mathqed{  \vrule width5pt height5pt depth0pt}

\def\V{W}

\def\df{d^{\phi}}
\def\hk{\_{\rm l}\,}
\def\n{\nabla}
\def\w{\wedge}

\def\cu{{\cal U}}   \def\cc{{\cal C}}   \def\cb{{\cal B}}  \def\cz{{\cal Z}}
\def\cv{{\cal V}}   \def\cp{{\cal P}}   \def\ca{{\cal A}}
\def\cw{{\cal W}}   \def\co{{\cal O}}
\def\ce{{\cal E}}   \def\ck{{\cal K}}
\def\ch{{\cal H}}   \def\cm{{\cal M}}
\def\cs{{\cal S}}   \def\cn{{\cal N}}
\def\cd{{\cal D}}
\def\cl{{\cal L}}
\def\cp{{\cal P}}
\def\cf{{\cal F}}
\def\ccr{{\cal  R}}

\def\gerG{{\fr{\hbox{g}}}}
\def\gerB{{\fr{\hbox{B}}}}
\def\gerR{{\fr{\hbox{R}}}}
\def\p#1{{\bf P}^{#1}}
\def\vf{\varphi}

\def\wt{\widetilde}
\def\wh{\widehat}

\def\and{\qquad {\rm and} \qquad}
\def\arr{\longrightarrow}
\def\ol{\overline}
\def\bbr{{\mathbb R}}\def\bbh{{\mathbb H}}\def\bbo{{\mathbb O}}
\def\bbc{{\mathbb C}}
\def\bbq{{\mathbb Q}}
\def\bbz{{\mathbb Z}}
\def\bbp{{\mathbb P}}
\def\bbd{{\mathbb D}}

\def\a{\alpha}
\def\b{\beta}
\def\d{\delta}
\def\e{\epsilon}
\def\f{\phi}
\def\g{\gamma}
\def\k{\kappa}
\def\l{\lambda}
\def\o{\omega}

\def\s{\sigma}
\def\x{\xi}
\def\z{\zeta}

\def\D{\Delta}
\def\L{\Lambda}
\def\G{\Gamma}
\def\O{\Omega}

\def\bd{\partial}
\def\bdf{\partial_{\f}}
\def\lag{Lagrangian}
\def\psh{plurisubharmonic }
\def\ph{pluriharmonic }
\def\pph{partially pluriharmonic }
\def\omp{$\omega$-plurisubharmonic \ }
\def\ffl{$\f$-flat}
\def\PH#1{\widehat {#1}}
\def\lloc{L^1_{\rm loc}}
\def\dbar{\ol{\partial}}
\def\lp{\Lambda_+(\f)}
\def\lpp{\Lambda^+(\f)}
\def\bo{\partial \Omega}
\def\Ob{\overline{\O}}
\def\fc{$\phi$-convex }
\def\PSH{{ \rm PSH}}
\def\SH{{\rm SH}}
\def\totr{ $\phi$-free }
\def\BM{\lambda}
\def\Der{D}
\def\CH{{\cal H}}
\def\RH{\overline{\ch}^\f }
\def\pconv{$p$-convex}
\def\MA{MA}
\def\lagpsh{Lagrangian plurisubharmonic}
\def\hermsk{{\rm Herm}_{\rm skew}}
\def\PSHl{\PSH_{\rm Lag}}
 \def\ppsh{$\pp$-plurisubharmonic}
\def\fp{$\pp$-plurisubharmonic }
\def\fh{$\pp$-pluriharmonic }
\def\Symn{\Sym(\rn)}
 \def\ci{C^{\infty}}
\def\USC{{\rm USC}}
\def\fa{{\rm\ \  for\ all\ }}
\def\ppc{$\pp$-convex}
\def\cpt{\wt{\cp}}
\def\ft{\wt F}
\def\ob{\overline{\O}}
\def\Be{B_\e}
\def\K{{\rm K}}

\def\M{{\bf M}}
\def\N#1{C_{#1}}
\def\ds{Dirichlet set }
\def\dir{Dirichlet }
\def\Fa{{\oa F}}
\def\TR{{\cal T}}
 \def\LAG{{\rm LAG}}
 \def\ISO{{\rm ISO_p}}
 \def\Span{{\rm Span}}

\def\AA{1}
\def\BB{2}
\def\CC{3}
\def\DD{4}
\def\EE{5}
\def\FF{6}
\def\GGG{7}
\def\HH{8}
\def\II{9}
\def\JJ{10}
\def\KK{11}
\def\LL{12}
\def\MM{13}

\vskip .4in

\def\E{E}
\def\bL{{\bf \Lambda}}
\font\headfont=cmr10 at 14 pt

\vskip .1in

%\begin{abstract}\end{abstract}

\title[REMOVABLE SINGULARITIES  FOR NONLINEAR SUBEQUATIONS]{REMOVABLE SINGULARITIES  \\ FOR NONLINEAR SUBEQUATIONS
}
 
\date{\today}
\author{ F. Reese Harvey and H. Blaine Lawson, Jr.$^1$}

\maketitle

%{\small\tableofcontents}

%\section{Introduction}
%\label{intro}

\centerline{\bf Abstract} \medskip
  \font\abstractfont=cmr10 at 10 pt
Let $F$ be a fully nonlinear second-order partial differential subequation
of degenerate elliptic type on a manifold $X$.  We study the question:
Which  closed subsets $E\ss X$ have the property that every $F$-subharmonic
function (subsolution) on $X-E$, which is locally bounded across $E$, extends to an $F$-subharmonic
function on $X$. We also study the related question for $F$-harmonic functions (solutions) which are continuous across $E$.
%Which  closed subsets $E\ss X$ have the property that every continuous function which is $F$-harmonic
% on $X-E$ extends to an $F$-harmonic function on $X$. 
The main result asserts that if  there exists  a convex cone
 subequation $M$ such that $F+M\ss F$, then any closed set $E$ which is $M$-polar has these properties.
$M$-{\sl polar}  means that $E=\{\psi=-\infty\}$ where $\psi$ is $M$-subharmonic on $X$ and smooth
 outside of $E$. Many examples and generalizations are given.  These include removable 
 singularity results for {\sl all branches} of the complex and quaternionic Monge-Amp\`ere equations,
and a general removable singularity result 
 for the harmonics of geometrically defined subequations.
 
 For pure second-order subequations in $\rn$ with monotonicity cone
 $M$ the {\sl Riesz characteristic} $p=p_M$ is introduced, and 
 extension theorems are proved for any closed singular set 
 $E$ of locally finite Hausdorff $(p-2)$-measure. This applies for example to branches of the
 equation $\s_k(D^2 u)=0$ ($k^{\rm th}$ elementary function) where $p_M=n/k$,
 and its complex and quaternionic counterparts where $p_M = {2n\over k}$ and $p_M={4n\over k}$ 
 respectively.
 
 For convex cone subequations themselves, several  removable singularity theorems are
 proved, independent  of the results above.

\vglue .9cm
\smallbreak\footnote{}{ $ {} \sp{ }{\rm Partially}$  supported by
the N.S.F. }

\vfill\eject

\ \ 
\vskip.6in

%{\small\tableofcontents}
\centerline{\bf TABLE OF CONTENTS} \bigskip

{{\parindent= .1in\narrower \noindent

\qquad \AA.     Introduction.    \smallskip

\qquad \BB.     Subequations and Subharmonic Functions.  \smallskip

\qquad \CC.     Removable Sets.  \smallskip

\qquad \DD.   Monotonicity Subequations.\smallskip

\qquad \EE.   Polar Sets for a Subequation.\smallskip

\qquad \FF.   Polar Sets are Removable.\smallskip

\qquad \GGG.   Harmonic Functions.\smallskip

\qquad \HH.  Classical Plurisubharmonic Functions.  

\smallskip

\qquad \II.   $p$-Plurisubharmonic Functions and $p$-Monotonicity.  

\smallskip

\qquad \JJ.  $p$-Pluripolar Sets and Riesz Potentials.

\smallskip

\qquad \KK. Subequations which are   $\cp_p$-monotone .

\smallskip

\qquad \LL.  Removable Singularities for Convex Subequations.

}}

{{\parindent= .9in\narrower

\bigskip
\noindent
Appendix A.\    A removable Singularity Theorem 
 of Caffarelli, Li and Nirenberg.
\smallskip

\noindent
Appendix B.\    An Illustration: Cones with $p={3\over2}$ and $n=3$.
\smallskip

}}

\vfill\eject

  %%%%%%%%%%%%%%%%%%%%%%%%%%%%%%%%%%%%%%%%%%%%%%%%%%%%%%%
  %%%%%%%%%%%%%%%%%%%%%%%%%%%%%%%%%%%%%%%%%%%%%%%%%%%%%%%
  %%%%%%%%%%%%%%%%%%%%%%%%%%%%%%%%%%%%%%%%%%%%%%%%%%%%%%%
  %%%%%%%%%%%%%%%%%%%%%%%%%%%%%%%%%%%%%%%%%%%%%%%%%%%%%%%
  %%%%%%%%%%%%%%%%%%%%%%%%%%%%%%%%%%%%%%%%%%%%%%%%%%%%%%%
  %%%%%%%%%%%%%%%%%%%%%%%%%%%%%%%%%%%%%%%%%%%%%%%%%%%%%%%

\noindent{\headfont \AA.\  Introduction}
\medskip
 
\noindent{\bf Some Historical  Background (Linear Equations).} Removable singularity results for solutions to {\sl linear} partial differential equations
 have a long history with a cornucopia of results.  The general set-up involves two
 things: first a closed subset $E$ of a manifold  $X$ which, since most results are local,
 might as well be an open subset of $\rn$, and second, a solution $u$ of the PDE in the complement
 $X-E$.  The conclusion is generally the same, namely, that  $u$ extends to a solution on $X$.
 However, two different  types of hypotheses are required, one on the ``growth'' of $u$ across $E$, and the 
 other on the ``size'' of the exceptional set $E$.
 
 The first classical result of this kind was the Riemann Removable Singularity Theorem
 for holomorphic functions on the punctured disk satisfying $f(z) = o(|z|^{-1})$.
 In 1956 Bochner generalized Riemann's result to solutions of {\sl  linear}  PDE's  which satisfy  the 
 ``growth hypothesis''  $f(x) =o(\dist(x,E)^{-q})$.  The remarkable thing about Bochner's 
 removable singularity theorem is that the only information needed about the differential equation
 is its order $m$.  The ``size'' hypothesis on $E$ is that in dimension $n-m-q$ the set 
 $E$ should be (in some sense) locally finite.  The majority of removable singularity
 results are of this kind -- only depending on the order $m$ of the equation.
 (See [\HP], [\HPP] for more detailed results.)

The main exception is where the ``growth'' assumption on $u$ states that $u$ is locally bounded
across $E$.  If $E=\{0\}$ in $\bbr^2$, then the harmonic functions extend, while the fundamental
solution for the wave equation, namely the characteristic function of the forward light cone,
does not extend, even though it is bounded.  In the bounded harmonic case, capacity  is the right notion
of ``size'', with capacity zero  for $E$  giving a necessary and sufficient condition for removability.

\medskip
The purpose of this paper is to prove removable singularity results for fully 
non-linear second-order equations which
are ``elliptic'' but can be very degenerate.  We begin with a removable singularity result
of a very general nature for subharmonics (or subsolutions) which are locally bounded across
the singular set $E$.   This theorem is then applied to yield the appropriate result, again of 
a very general nature, for  the harmonics  (or viscosity solutions). These results are then
applied -- worked out in detail and sometimes  generalized -- in a number of special cases,  including all the branches of the
complex and quaternionic Monge-Amp\`ere equations, and also all the geometrically defined subequations.

We also obtain extension results across closed sets of locally finite Hausdorff $k$-measure for computable
values of $k$. Here again there are many interesting applications.
These are the first such results for general non-uniformly elliptic equations.

The central concept here is that of a subequation $F$ (see Section \BB)
{\bf  which is $M$-monotone} for some convex cone
subequation $M$, that is, $F+M = F$ point-wise on $X$. A  subset $E\ss X$ is called
{\sl $C^\infty$ $M$-polar}  if $E=\{\psi=-\infty\}$ for an $M$-subharmonic function $\psi$ which is
smooth in $X-E$. The first main result is the following.

\Theorem{\FF.1} {\sl
Suppose $F$ is a subequation on $X$ which is $M$-monotone.   If $\E\ss X$ is a closed subset 
which is locally $C^\infty$ $M$-polar,
then $\E$ is removable for $F$-subharmonic functions which are locally bounded above across $\E$.
}

\medskip
The hypothesis that  the function be locally bounded above across $E$ is automatic in
certain cases (see Lemma \FF.2).

Theorem \FF.1 gives a corresponding result for $F$-harmonics.
 These are functions $u$ such that $u$ is a subsolution (subharmonic) for $F$ 
 and $-u$ is a subsolution  for the dual subequation $\ft$ (see (\GGG.1)).
 
 \Theorem {\GGG.1}  {\sl
Suppose $F$ is a subequation which is $M$-monotone
and $\E$ is a closed set with no interior, which  is locally  $C^\infty$  $M$-polar.
Then for $u \in C(X)$,

\medskip
\centerline{
$u$ is  $F$-harmonic   on $X-\E \qquad\Rightarrow\qquad u$ is $F$-harmonic   on $X$.
} 
}

 \medskip

 These two theorems have a number of  applications,  explored in Sections
\HH-\KK.
 
When $M$ is a convex cone subequation, one has $M+M\ss M$ and
 $\wt M+ M \ss \wt M$.  Applying Theorem \FF.1 with $F=M$ should be viewed as a 
 weak result. In general there are plenty of removable sets for $M$-subharmonics which
 are not $M$-polar (See Section \LL). However, applying Theorem \FF.1 with 
 $F=\wt M$ is a much better result since $\wt M \supset M$  is  much
 larger than $M$ and  highly non-convex. See Remark \GGG.2 and Corollary \GGG.3

 A basic and illuminating case is that of  classical plurisubharmonic functions in $\bbc^n$.
 Here   $M$ is  the convex cone subequation $\cp^\bbc$  defined by requiring 
  the complex hermitian hessian $(D^2 u)_\bbc \cong (\partial^2u/\partial z_i \partial \bar z_j)$ to be $\geq0$.
 The $\cp^\bbc$-polar sets are the classical pluripolar sets, and the theorems above apply
 to prove new removable singularity results for {\bf all branches} of the complex Monge-Amp\`ere equation (Theorem \HH.1). For the principal branch, the $\cp^\bbc$-harmonic functions are just
 the continuous maximal functions, i.e., the continuous plurisubharmonic functions $u$ satisfying
 $\det_\bbc(D^2 u) = 0$, and there are extendability results due to Bedford and Taylor  [\BTA], [\BTB], [\BTC].

 Another important case involves the $p$-plurisubharmonic functions, 
 which are discussed in Sections \II, \JJ\ and \KK.   When $p$ is an integer,
these are the  functions which become classically subharmonic when restricted to 
minimal $p$-dimensional submanifolds.
They are studied extensively in [\pPSH].
 The corresponding convex cone subequation $\cp_p$ is a monotonicity cone for a large collection
 of interesting subequations and their duals (see Section \KK).  
 For appropriate integers $p$ these include the Lagrangian and Special Lagrangian subequations
 in $\bbc^n$, the associative and coassociative subequations in $\bbr^7$,  the Cayley subequations
 in $\bbr^8$, and their analogues on riemannian manifolds (see [\DDR]).
 %In fact, the main results (Theorems \II.2 and \JJ.5) apply to all the
  %``geometric'' or  $\GG$-subequations yielding results
% for $\GG$-harmonics and dually $\GG$-plurisubharmonics. 
 The cases where $p$ is not an integer also have important applications, and lead
 to the concept of Riesz characteristics, discussed below.

 The $\cp_p$-polar sets are studied in Section \JJ \  using Riesz potentials, and the following is proved.
 
\medskip
\noindent
{\bf Theorem  \JJ.5.}
{\sl
Any closed subset $\E\ss \rn$ with locally finite  Hausdorff $(p-2)$-measure,  is locally $C^\infty$ 
$\cp_p$-polar.  
More generally, any closed subset $E\ss \rn$ with $(p-2)$-capacity zero is locally 
$C^\infty$ $\cp_p$-polar.
Therefore, if a subequation $F$ is $\cp_p$-monotone, then any such set
 is removable for $F$-subharmonics and $F$-harmonics, 
as in the theorems above.
}
\medskip

Any subequation $F=\cp(\GG)$ which is defined geometrically by a closed
subset $\GG\ss G(p,\rn)$ of the Grassmannian of $p$-planes (see (\KK.1))
 is $\cp_p$-monotone and the theorem applies. It gives extension theorems 
 for $\GG$-harmonics and $\GG$-plurisubharmonics, and dually $\GG$-plurisubharmonics. 
 This includes the specific geometric
 cases cited above. For example, for Lagrangian or Special Lagrangian subequations in $\bbc^n$ 
 closed sets of locally finite Hausdorff $n-2$-measure are removable.
 
 More generally let $F\ss\Symn$ be a pure second-order 
 subequation with monotonicity cone $M$.  If $\cp_p\ss M$, then $F$ is $\cp_p$-monotone, and
Theorem \JJ.5 applies. This leads to the question of determining the largest $\cp_p \ss M$.
An easily computable invariant, called the {\bf Riesz characteristic $p_M$ of $M$}, provides the answer
(see Definition \KK.3 and Theorem \KK.4).

Theorem \JJ.5 then gives the following (cf. Theorem \KK.6).

\medskip\noindent
{\bf Theorem A.}  {\sl
Suppose $F\ss\Symn$ is a subequation in $\rn$ with monotonicity cone $M\ss\Symn$.
 Let $p=p_M$ be the Riesz characteristic of $M$. 
 Then any closed subset $E$ of locally finite
Hausdorff $(p-2)$-measure (or, more generally, of $(p-2)$-capacity zero) is  removable for $F$-subharmonics and $F$-harmonics, 
as in the theorems above.
}
\medskip

A large set  of interesting subequations to which these results apply is
the family defined by homogeneous polynomials $P$ on $\Symn$ which
are G\aa rding hyperbolic with respect to the identity (cf. [\GAR],  [\HYP], [\HLGAR]).
Here the convex cone subequation $M$ corresponds to the G\aa rding cone,
and each branch of the equation $P(D^2 u)=0$ is $M$-monotone.

More specifically, a homogeneous polynomial $P:\Symn \to \bbr$ of degree $m$  is
 {\bf G\aa rding hyperbolic} with respect to $I$ if $P(I)>0$ and 
 for each $A$, all the roots of the one-variable polynomial $p_A(t) = P(A+tI)$ are real.
 In this case one can write $p_A(t) = \prod_i(t+\l_i^P(A))$ where $\l_1^P(A) \leq \cdots \leq  \l_m^P(A)$ 
are the {\sl ordered $P$-eigenvalues} of $A$. We assume $\l_i^P(A+B)\geq \l_i^P(A)$ whenever
$B\geq0$. Then each {\sl branch}
$$
\L^P_i \ \equiv\ \{ A : \l^P_i(A)\ \geq\ 0\}, \qquad 1\ \leq\  i\ \leq \ m
$$
is a subequation, and the principal branch
$$
M  \ \equiv \ \L^P_1 \ =\ \{ A : \l^P_i(A)\ \geq\ 0\ \forall\, i\},
$$
called the {\bf G\aa rding cone}, is a {\bf convex} cone which satisfies $\L^P_i + M\ss \L^P_i$,
 that is, it is a monotonicity cone for  each branch of the equation $P(D^2 u)=0$.

A basic example is given by the $k^{\rm th}$ elementary symmetric function $\s_k(A)$
of the eigenvalues of $A$.  The G\aa rding cone is 
$M= \Sigma_k = \{A : \s_1(A)\geq0, ..., \s_k(A)\geq 0\}$, and the Riesz characteristic is computed to be 
$$
p_{\Sigma_k} \ =\ {n\over k}.
$$
Note that  $\Sigma_1=\cp$,  the basic homogeneous  real Monge-Amp\`ere equation
on convex functions, whose study goes back to Alexandrov [\ALEX] and Pogorelov [\POGA], [\POGB].
The other G\aa rding cones $\Sigma_k$, sometimes called {\sl Hessian equations},
have been studied by Trudinger-Wang [\TWA], [\TWB], [\TWC],  Labutin [\LAB], and others. 
Note however, that each $\Sigma_k$ has $k-1$ other  branches,
which are neither convex nor uniformly elliptic, but  to which 
Theorem A applies with $p=n/k$.

Further interesting G\aa rding polynomials are discussed in Section \KK.
One is the $p$-fold sum operator, which for integer
$p$ is given by
$
M_p(A) = \prod ( \l_{i_1}(A)+\cdots + \l_{i_p}(A))
$,
whose principal branch  is $\cp_p$.

We note that every subequation $F\ss\Symn$, which is defined purely in terms
of the eigenvalues of the matrices (i.e.,  is O($n$)-invariant), has {\bf complex and quaternionic
analogues}  $F^\bbc \ss \Sym_\bbr(\bbc^n)$ and $F^\bbh \ss \Sym_\bbr(\bbh^n)$.
The set $F^\bbc$ is defined by imposing the constraints of $F$ on the $n$ eigenvalues of the hermitian
symmetric part $A_\bbc = \half(A-JAJ)$ of $A$. The set $F^\bbh$ is defined similarly.
This applies for example to the equations $\Sigma_k$ above, and the special cases
$\cp^\bbc=\Sigma^\bbc_1$ and $\cp^\bbh=\Sigma^\bbh_1$ correspond to the complex and quaternionic Monge-Amp\`ere
equations on the underlying  plurisubharmonic functions.  In complete generality the 
Riesz characteristics satisfy:
$$
p_{F^\bbc} \ =\ 2 p_F  \and 
p_{F^\bbh} \ =\ 4 p_F
$$

\medskip

\noindent
{\bf Note (Uniform ellipticity).} There are many natural cones, which when they
are monotonicity cones for a subequation $F$, imply that $F$ is uniformly 
elliptic. One family of such cones is  $\cp(\d) = \{A :  A+\d (\tr A) I\geq0\}$. Another is the
family of 
{\sl Pucci cones} $\cp_{\l,\L}$ defined in Example 2 of Section \KK.  
Monotonicity with respect to $\cp_{\l,\L}$
corresponds to the standard $(\l,\L)$-uniform ellipticity.
For $\cp(\d)$ the Riesz characteristic is $(1+\d n)/(1+\d)$, and for $\cp_{\l,\L}$
the Riesz characteristic is $p_{\l,\L}={\l\over\L} (n-1)+1$.  Thus as the uniform ellipticity constant
increases, the Hausdorff dimension of the removable sets increases.
However, we point out that the main thrust of this paper is to establish 
removable singularity results for equations which are {\sl not} necessarily uniformly elliptic.
\medskip

\noindent
{\bf Further Results -- Variable Coefficients.} There  is a theorem similar to Theorem A  above which holds for variable-coefficients (see the end of Section \KK).  There is also a result on general
riemannian manifolds where $\cp_p$ is replaced by its riemannian analogue,
defined using the riemannian hessian.  See Theorem \KK.7.

\medskip

After circulating an early version of this paper, we received an article [\AGV] 
by Amendola, Galise and Vitolo which, in the language above, studied 
uniformly elliptic equations of Riesz characteristic $p_{\l,\L}\geq2$
(see condition (1.1) in [\AGV]). The authors proved the extendability of 
the maximum principle and comparison, and also the removability of 
singularities across sets of $(p_{\l,\L}-2)$-capacity zero.
In particular, they established Theorem  A above in the $(\l,\L)$-uniformly
elliptic case.

Removable singularity theorems have also been established
for certain uniformly elliptic equations by Labutin (see [\LAB], [\LLAB], [\LLLAB]
and references therein.)
 Results on the  removability of isolated singularities for certain uniformly 
 elliptic equations also appeared in papers of Armstrong, Sirakov and Smart [\ASSA], [\ASSB].
There are also some very recent results of Vitolo [\VIT].
\medskip

For convex cone subequations $M$ themselves, a wide variety of strong removability
results are obtained rather easily from the linear case 
using the Strong Bellman Principle [\ASPECTS], provided that
$M$ is ``second-order complete''. Using linear results from [\HP] (see also [\HAR], [\HPP]),
some of these possiblities are discussed in Section \LL.
None of them depend on the rest of the paper.
For example, we show that every closed set $E$ of locally finite
Hausdorff codimension-2 measure is removable for $M$-subharmonic functions
which are locally bounded above across $E$.  Several other removable singularity results,
which entail various combinations  of ``growth'' and ``size'', are also given 
 in Section \LL.

Results in this paper were inspired by work of Caffarelli, Li and Nirenberg [\CLN].
In Appendix A we present an alternate proof of one of their results.

%%%%%%%%%%%%%%%%%%%%%%%%%%%%%%%%%%%%%%%%%%%%%%%%%%
%%%%%%%%%%%%%%%%%%%%%%%%%%%%%%%%%%%%%%%%%%%%%%%%%%
%%%%%%%%%%%%%%%%%%%%%%%%%%%%%%%%%%%%%%%%%%%%%%%%%%
%%%%%%%%%%%%%%%%%%%%%%%%%%%%%%%%%%%%%%%%%%%%%%%%%%
%%%%%%%%%%%%%%%%%%%%%%%%%%%%%%%%%%%%%%%%%%%%%%%%%%

 \vfill\eject
% \vskip .3in

\noindent{\headfont \BB.     Subequations and Subharmonic Functions. }
\medskip
 
By a  {\bf subequation} on a manifold $X$ we mean  a closed subset of the 2-jet bundle of $X$,
(locally, a closed subset of $X\times \bbr\times\rn\times \Symn$ with
$X^{\rm open}\ss \rn$) with non-empty fibres  over $X$ satisfying the primary  condition
$$
F+\cp\ \ss\  F
\eqno{\bf (Positivity)}
$$
where $\cp = \{(x,0,0,A) : x\in X\ {\rm and}\ A\geq0\}$.
Intrinsically, $\cp$ is the set whose fibre at $x$ consists of 2-jets of functions with minimum value
0 at $x$.
In [\DDR, \S\S2 and 3] a subequation is  assumed to satisfy  a  further 
 mild topological regularity  condition which  is important for discussing the dual subequation $\wt F$ (see Section \GGG).
    If each fibre $F_x$ of $F$ is a convex cone with vertex at the origin
in $\bbr\times\rn\times\Symn$,
then $F$ is called a {\sl convex cone subequation}.

Let $\USC(X)$ denote the space of upper semi-continuous $[-\infty,\infty)$-valued functions $u$
on $X$. Often $u$ is not allowed to be $\equiv -\infty$ on any connected component of $X$.
 By a   {\bf   (viscosity)   test function} for $u\in \USC(X)$ at a point $x$  we mean a
$C^2$-function $\vf$ defined  near $x$ such that  $u-\vf$ has local maximum value 0 at $x$.
(See [\CIL], [\CRA] and   references therein for a discussion of the viscosity approach to differential equations.)

\Def{\BB.1} A function $u\in \USC(X)$ is {\bf  $F$-subharmonic on $X$} if for each
$x\in X$ and each test function $\vf$ for $u$ at $x$, we have
$$
(\vf(x), D_x\vf, D_x^2\vf) \ \in\ F_x.
$$

  Positivity  ensures that if $u\in C^2(X)$ and
$(u(x), D_x u, D_x^2 u)\in F_x$ for all $x$, then $u$ is $F$-subharmonic on $X$.
Each subequation $F$ has its own potential theory, which studies the set $F(X)$ of 
$F$-subharmonic functions. Of course, each such potential theory has its own unique characteristics.
However, there is a surprisingly large common core of general results.  See, for example,
[\DDR, Theorem 2.6] for the beginning of a list of results that hold for all these theories.
The purpose of this note is to add a removable singularity theorem to each of these  potential theories.

\vskip .3in

\noindent{\headfont \CC.     Removable Sets. }

\medskip
We consider closed subsets $\E$ of $X$, and we restrict 
 attention to  upper semi-continuous functions $u$
 on $X-\E$ which are locally bounded above at points of $\E$.   Each such function $u$
 has a canonical upper semi-continuous extension $U$ across $\E$ to all of $X$ defined as 
 follows.  If $\Int \E = \emptyset$, set
 $$
 U(x) \ \equiv \ \limsup_{ 
  \begin{matrix} 
  y\to x  \cr  y\notin \E
  \end{matrix}} u(y)
 \ =\ \lim_{\e\searrow 0} \sup_{B(x,\e)-\E} u
 \eqno{(\CC.1)}
 $$
When $\Int \E \neq  \emptyset$, extend this by setting $U(x)\equiv -\infty$ on $\Int \E$.
It is easy to see that $U$ is the upper semi-continuous regularization of the function 
$v$, which is defined to be $u$ on $X-\E$ and $-\infty$ on $\E$.
(This fact will be used in Section \FF \ to complete the proof of Theorem \FF.1.)

\Def{\CC.2} A closed  set $\E$ is {\bf removable for $F$-subharmonic functions 
locally bounded above across $\E$} if for each function $u$
which is $F$-subharmonic on $X-\E$ and  is  bounded above in a neighborhood of each point of  $\E$,
the canonical extension $U$ is $F$-subharmonic on $X$.
\medskip

The most elementary removable singularity result can be stated as follows.

\Prop{\CC.3} {\sl
If $U\in\USC(X)$ and $U$ has no test functions at points $x\in E$, then 
$U\in F(X-E)\ \ \Rightarrow\ \ U\in F(X)$.
}
\medskip

For example, take $U(x) = u(x) +\l|x|$ with $u$ smooth and $E=\{0\}$.

\vskip .3in

\noindent{\headfont \DD.   Monotonicity Subequations.}
\medskip

Given a subequation $F$ on $X$, a {\bf monotonicity subequation} for $F$ 
is a convex cone subequation $M$ on $X$ satisfying
$$
F+M\ \ss\ F
\eqno{\bf (Monotonicity)}
$$
where ``+'' denotes fibre-wise sum.  In this case we say the subequation $F$ is {\bf $M$-monotone}.
It follows easily from the definition of $F$-subharmonicity that:
$$
u\in F(X)\ \ {\rm and}\ \ \psi\in M(X) \cap C^\infty(X)     \qquad\Rightarrow\qquad 
u+\e\psi \in F(X) \ \ {\rm for\  all}\ \ \e>0.
\eqno{(\DD.1)}
$$

\vskip .3in
%\vfill\eject

\noindent{\headfont \EE.   Polar Sets for a Subequation.}
\medskip
 
For the purposes of dealing with removable singularities
it is convenient to limit attention to a restricted class of ``polar sets''  for a subequation.

A function $\psi \in\USC(X)$  is called  a {\bf polar function} for a closed set  $\E$  if 
$\E = \{x\in X:\psi(x) = -\infty\}$.  If in addition $\psi$ is
smooth on  $X-\E$, it is called a {\bf smooth polar function} for $\E$.
Let $M$ be a subequation.

\Def{\EE.1}  A closed subset $\E \ss X$ is called $C^\infty$ {\bf $M$-polar} if it  admits a smooth
polar function which is $M$-subharmonic on $X$.  

\medskip

%In the following discussion this concept will only be utilized when $F$ is a convex cone
%subequation such as the $M$ above.

%This concept will only be utilized when $F$ is a convex cone
%subequation such as the $M$ above.

%We shall consider such sets only for convex cone subequations such as the $M$ above.

%\vskip .3in
\vfill\eject

\noindent{\headfont \FF.   Polar Sets are Removable.}
\medskip
 
This is a classical result for the euclidean Laplacian.  Moreover,  the classical proof
extends quite easily to the following very general case. Let $M$ be a convex cone subequation.

\Theorem{\FF.1} 
{\sl
Suppose $F$ is a subequation which is $M$-monotone.  If $\E$ is a closed subset 
which is locally $C^\infty$ $M$-polar,
then $\E$ is removable for $F$-subharmonic functions which are locally bounded above across $\E$.
}

\pf
Given $u\in F(X-\E)$ we must show that  $U\in F(X)$.  Let $\psi$ denote a smooth $M$-subharmonic
polar function for $\E$.  As noted above in (\DD.1), each $U+\e\psi$ is $F$-subharmonic on $X-\E$.
Since $(U+\e\psi)(x) = -\infty$ at each point $x\in \E$, $U+\e\psi$ cannot have a test function
at such a point $x$. This proves that for each $\e>0$, $U+\e\psi \in F(X)$.  Now the set
$\{\psi<0\}$ is an open set containing $\E$, and so we can assume that $X=\{\psi<0\}$, i.e., 
$\psi <0$ on $X$.

By the ``families bounded above property'' for $F$-subharmonic functions ([\DDR, Thm. 2.6 (E)]),
the upper semi-continuous regularization of the upper envelope $v$
of the family $\{U+\e\psi\}$ is $F$-subharmonic on $X$. This upper envelope $v$ is given by
$$
v(x)\ =\ \begin{cases} u(x) \qquad {\rm if}\ \ x\notin \E  \cr
-\infty \qquad {\rm if}\ \ \ x\in \E  
\end{cases}
$$
since $\E$ is exactly the $-\infty$ set of $\psi$.
Finally, as noted in Section \CC,  $v$ has upper semi-continuous regularization $U$.\qed

\medskip

The hypothesis in Theorem \FF.1 of being bounded above across $E$ can be deleted for
certain subequations.  Applications of this result will be left to the reader.

\Lemma{\FF.2}  {\sl
Suppose $F\ss \L_p \equiv \{\l_p(A) \geq0\}$ and $E$ has $(n-p)$-dimensional Hausdorff
measure zero.  Then each $F$-subharmonic function on $X-E$ is locally bounded above across
$E$.
}

\pf
Pick a point in $E$ which we can assume is the origin.  As in Shiffman [\SHF]
we can construct a product of balls $B\times B'\ss\rn$
with $B$ a ball in $\bbr^p$ such that $E\cap (\partial B\times B')=\emptyset$
and the projection of $E\cap(B\times B')$ on $B'$ has measure zero.  
By the Restriction Theorem 5.3 in [\REST], the function 
 $u$ restricted to the good slices $B\times \{y\}\ss B\times B'-E$ 
 is a subaffine function, and hence satisfies the maximum principle on
 each good slice.  Therefore, 
 $$
 u(x,y) \leq \sup_{x\in\partial B} u(x,y) \qquad {\rm for\ almost\ every\ } y,  \ \ {\rm and\ so}
 $$
$$
 \sup_{B\times B' - E} u \   \leq \sup_{\partial B\times B'} u \ <\ \infty. \qquad
 \mathqed
 $$

\vskip .5in
%\vfill\eject

\noindent{\headfont \GGG.   Harmonic Functions}
\medskip

Given a subequation $F$ recall the {\bf dual subequation}
$$
\wt F \ =\  -(\sim \Int F) \ =\  \sim(-\Int F)
\eqno{(\GGG.1)}
$$
 (see [\DDD], [\DDR]).  
Here a mild topological condition (T) on the set $F$ is required [\DDR, 3.2]. 

\medskip
\noindent
{\bf Definition.}
A function $u$ is
{\bf $F$-harmonic} on $X$ if $u\in F(X)$ and $-u\in \wt F(X)$.  
\medskip

For smooth functions
$u$ this means that 
$$
(u(x), D_xu, D^2_xu)\in \partial F_x \qquad {\rm  for\  each  \ } x\in X.
$$

Let $M$ be a convex cone subequation.

\Theorem {\GGG.1}  {\sl
Suppose $F$ is a subequation which is $M$-monotone
and $\E$ is a closed set with no interior, which  is locally  $C^\infty$  $M$-polar.
Then for $u \in C(X)$,

\medskip
\centerline{
$u$ is  $F$-harmonic   on $X-\E \qquad\Rightarrow\qquad u$ is $F$-harmonic   on $X$.
} 
}

\pf
Corollary 3.5 in [\DDR] states that
$$
F+M\ \ss\ F \qquad\iff\qquad \wt F + M \ \ss\ \wt F
\eqno{(\GGG.2)}
$$
for any subset $M$ of the 2-jet bundle.   Theorem \FF.1 applies 
to the $F$-subharmonic function $u$ on $X-\E$
proving that the canonical extension $U$ defined by (\CC.1) is $F$-subharmonic on $X$.
Similarly, because of (\GGG.2),   the $\wt F$-subharmonic function
$v=-u$ on $X-\E$ has canonical extension $V$ which is $\wt F$-subharmonic on $X$.
  By continuity of $u$ on $X$  and the fact that $\Int \E=\emptyset$, the 
canonical extensions $-U$ and $V$ agree,
proving that $U$ is $F$-harmonic on $X$.  \qed

 \Remark{\GGG.2}  Note that for a convex cone subequation $M$ the space of 
 dual $\wt M$-subharmonics is much larger than the space of  $M$-harmonics.
 To see this note that (\GGG.2) implies that $\wt M$ is $M$-monotone, and since 
 $0\in \wt M$, it follows that  $M\ss\wt M$.  Because of this, the special case
  of Theorem  \FF.1, stating that $M$-polar sets
 are removable for $M$-subharmonics, should be considered a weak result.
 In general there are plenty of removable sets for $M$-subharmonics which
 are not $M$-polar (See Section \LL).
 By contrast,  that $M$-polar sets are removable 
 for $\wt M$-subharmonics (Theorem \FF.1) and for $M$-harmonics
 (Theorem \GGG.1) constitute much better results, which are stated as follows.
 
 \Cor{\GGG.3}
 {\sl
 Suppose that $M$ is a convex cone subequation and $E$ is a closed set which
 is locally  $C^\infty$ $M$-polar.  Then $E$ is removable for $\wt M$-subharmonic functions which are 
 locally bounded above across $E$.  Moreover, if $E$ has no interior and $u\in C(X)$, then
 
 \medskip
 \centerline
 {
 $u$ is $M$-harmonic on $X-\E\qquad\Rightarrow\qquad 
 u$ is $M$-harmonic on $X$.
 }
 }

\medskip

\Remark{\GGG.4. (Horizontally Subharmonic)}   This example illustrates the general nature of 
Theorems \FF.1 and \GGG.1 and may help illuminate some of the trivialities that are involved.
Suppose $\rn = \bbr^p\times \bbr^q$.  Define $F$ by requiring that the horizontal Laplacian
 satisfy $\sum_{k=1}^p {\partial^2 u\over \partial x_k^2} \geq0$.
 Assume $q>0$, i.e., $F$ does not involve all the   variables in $\rn$.
Choose an open set $Y\ss \bbr^q$ and consider an arbitrary function $u: Y\to [-\infty,\infty)$.
Then:

\smallskip
\centerline
{
$u(y)$ is $F$-subharmonic on $\bbr^p\times Y \qquad\iff\qquad u\in \USC(Y)$
}
\smallskip
\centerline
{
$u(y)$ is $F$-harmonic on $\bbr^p\times Y \qquad\iff\qquad u\in  C(Y)$
}
\medskip
Note that this example is self-dual, i.e., $\wt F=F$ (if $p\geq1$)
and a convex cone subequation. 
It is now easy to see that the removable singularities result Theorem \FF.1 must apply to general 
functions $u\in \USC(Y-\E)$ yielding $U\in \USC(Y)$, while the hypothesis
in Theorem \GGG.1, in addition to $u\in C(Y-\E)$ must include enough
information to conclude that $U\in C(Y)$. 

Consider, for example,  
 $u(y) =\sin(1/|y|)$ with $\E=\{y=0\}$.
This function is $F$-harmonic on $\rn-\E$ and bounded across $\E$ but cannot be extended
to an $F$-harmonic function on $\rn$.  It is also easy to see that for this 
convex cone subequation, the hypothesis ``$\Int\E=\emptyset$'' occuring  in the second
half of Corollary \GGG.3  cannot be dropped.

\Remark{\GGG.5.(Convex Cone Subequations  which are Second-Order Complete)} 
If $M$ is such a subequation 
(see Section \LL), the $M$-subharmonic functions are also $L$-subharmonic for a linear equation
$L$ with positive definite second-order part.  This is enough to ensure that such a function is locally 
Lebesgue integrable.  In particular, this excludes the possibility  of a closed set $E$, 
which admits a smooth $M$-polar function, having any interior. 
Hence,  in this case the hypothesis $\Int \E = \emptyset$ can be dropped 
in Theorem \GGG.1 and Corollary \GGG.3.

\medskip

Here is a situation where the hypothesis of continuity on $u$ in Theorem \GGG.1
can be weakened to local boundedness.

\Theorem{\GGG.6} {\sl
Let $F$ be a subequation on $X$ which is $M$-monotone.
Let  $E\ss X$ be a compact $M$-polar set with no interior. 
Suppose that $E\ss \O^{\rm open}\ss\ss X$ where
the Dirichlet problem for $F$-harmonic functions on $\O$ is uniquely 
solvable for all continuous boundary functions.  Then if $u\in \USC(\overline \O -E)$ is $F$-harmonic on $\O-E$ and locally 
bounded across $E$, then $u$ extends to an $F$-harmonic function on $\O$.
}

\pf
Let $U=u^*$ be the upper semi-continuous regularization of $u$ on $\O$.
Let $-V = (-u)^* = -u_*$ be the  upper semi-continuous regularization of $-u$ on $\O$.
Then $V\leq U$, and by Theorem \GGG.1 we have that

\medskip
\centerline
{
$U$ is $F$-subharmonic on $\O$ \ \ \ and\ \ \  $V$ is $\ft$-subharmonic on $\O$.
}
\medskip
\noindent
Let $\wt u$ be the unique $F$-harmonic extension of $u\bigr|_{\bo}$  to $\O$.  Then
since $\wt u$ and $-\wt u$ are obtained by the Perron process, we have
$U\leq \wt u$ and $-V\leq -\wt u$. Hence, $\wt u\leq V\leq U \leq \wt u$.\qed

%\vfill\eject
\vskip.3in

\noindent{\headfont \HH.   Classical Plurisubharmonic Functions.}
\medskip

 The classical case is that of  plurisubharmonic functions and  pluripolar sets in $\bbc^n$.
  Note that in our parlance the plurisubharmonic functions are defined by the 
 constant coefficient, pure second-order subequation
  $$
 \cp^\bbc \ \equiv\ \{A\in\Sym_\bbr(\bbc^n) : A_\bbc \geq0\}
 $$
 where $A_\bbc \equiv \half(A-JAJ)$ denotes the hermitian symmetric component of $A$.
 (Here $J:\bbr^{2n}\to\bbr^{2n}$ denotes the complex structure.)
 Thus $u\in\USC(X)$ is plurisubharmonic if $D^2_x\vf   \in  \cp^\bbc$
 for each test function $\vf$ for $u$ at points $x\in X$.
 Pluripolar sets are classically  defined simply as subsets of  $-\infty$ level sets of  plurisubharmonic functions
 (not allowed to be identically $-\infty$ on any component of $X$). 
The pluripolar sets   are very well understood.

Fix an open set  $X\ss \bbc^n$ and consider the homogeneous complex Monge-Amp\`ere 
equation 
$$
\det_\bbc\left \{( D^2 u)_\bbc \right\}\ =\ 0
$$
Corresponding to this equation there are $n$ distinct, pure second-order subequations, or ``branches'', 
defined using the ordered eigenvalues $\l_1\leq\l_2\leq\cdots\leq\l_n$ of $A_\bbc  \equiv\ ( D^2 u)_\bbc$
 by setting
$$
\bL_{k}^\bbc \ \equiv\ \{ \l_k(A_\bbc)\geq0\}
$$
A $C^2$-function  $u$ which  is $\bL_{k}^\bbc$-harmonic for any one of these branches  satisfies the differential equation $\det_\bbc\{(D^2 u)_\bbc\} =0$, but there are $n$ distinct notions of being
$\bL_{k}^\bbc$-harmonic.  Only the real parts of holomorphic functions are $\bL_{k}^\bbc$-harmonic
for all $k$.

 Note that $\l_1(A_\bbc) \geq 0$ is the requirement that $A_\bbc \geq0$
so that $\cp^\bbc$ is the first, or smallest branch $\bL_{1}^\bbc$.
This subequation $\bL_{1}^\bbc=\cp^\bbc$ is the classical homogeneous  Monge-Amp\`ere 
equation for plurisubhamonic functions. 
It is the only convex branch.

One can show that a function is $\cp^\bbc$-harmonic if and only if it is continuous and maximally plurisubharmonic, ensuring that maximality is a local concept for continuous functions.

It is easy to see that the dual of the $k^{\rm th}$ branch $\bL_{k}^\bbc$ is $\bL_{n-k+1}^\bbc$.
In particular, the dual of the smallest branch $\cp^\bbc = \bL_{1}^\bbc$ is the largest
branch $\wt \cp^\bbc = \bL_{n}^\bbc$ which only requires that one eigenvalue $\l_{\rm max}$
 of $\left( D^2 u  \right)_\bbc$  be $\geq0$.
 (This largest branch $\wt \cp^\bbc$ is the complex analogue of the subequation $\cpt$
 whose subharmonics are the subaffine functions introduced in [\DDD].
 In fact, the $\cpt^\bbc$ -subharmonics are characterized as being ``sub'' the real parts of holomorphic
 functions on $\bbc^n$ [\REST,  Prop. 5.4].)
One can easily show that $\cp^\bbc$ is a monotonicity
cone for each subequation $\bL_{k}^\bbc$.  
As noted in Remark \GGG.5,  plursubharmonic functions are locally Lebesgue integrable.
Thus,  pluripolar sets have measure 
zero and hence no interior.
%Therefore, as a special case of Theorems \FF.1 and \GGG.1,
 %we have the following.

\Theorem{\HH.1} 
{\sl
Let $\E^{\rm closed}\ss X^{\rm open} \ss\bbc^n$ be a locally pluripolar set. 
Any function $u\in \bL_k^\bbc(X-\E)$
which is locally bounded above across $\E$ has a canonical extension  $U\in\bL_k^\bbc(X)$.
Moreover, if a   $\bL_k^\bbc$-harmonic function on $X-\E$ has a continuous extension to $X$,
then this extension is $\bL_k^\bbc$-harmonic   on $X$.
}

\medskip
\noindent
{\bf Proof of Theorem \HH.1}  If $\E$ is contained in $\{\psi = -\infty\}$ (locally) where $\psi$ plurisubharmonic
on $X$ and  smooth outside $\{\psi = -\infty\}$, and where $\{\psi = -\infty\}$ closed, then we can
simply replace $\E$ by $\{\psi = -\infty\}$ and apply Theorems \FF.1 and \GGG.1.
However, for the more general assertion of Theorem \HH.1 we proceed as follows.

Suppose that locally we have  a plurisubharmonic function $\psi$ with $\E\ss \{\psi = -\infty\}$
(but $\psi$ may not be smooth outside $\{\psi = -\infty\}$). The proof of Theorem \FF.1 is valid if we can show
that $u\in \bL_k^\bbc(X-\E)$ implies that $u+\e\psi \in  \bL_k^\bbc(X-\E)$.  This can be reduced to the 
case where $u$ and $\psi$ are quasi-convex (see, for example, Theorem 8.2 in [\DDD]).
In the quasi-convex case $u$ and $\psi$ have second derivatives $D_x^2u$ and $D_x^2\psi$ a.e.,
 and  $D_x^2(u+\e\psi) = D_x^2u + \e D_x^2\psi$ a.e.  By Corollary 7.5  in [\DDD] the fact that 
$D_x^2(u+\e\psi) \in \bL_k^\bbc +\cp^\bbc = \bL_k^\bbc$ a.e. implies that $u+\e\psi \in \
\bL_k^\bbc(X-\E)$.\qed

\vskip.3in

\noindent{\headfont \II.   $p$-Plurisubharmonic Functions and $p$-Monotonicity.}
\medskip

This is our second example, and again the roots are classical.
Fix  $p$ with  $1\leq p\leq n$.
First we consider the case where $p$ is an integer.
The  convex cone subequation
$\cp_p \equiv \cp(G(p,\rn))$ is defined in [\DDD] by
$$
\cp_p \ \equiv\ \{A\in\Symn : \tr_W A \geq0 \ \forall\, W\in G(p,\rn)\}.
\eqno{(\II.1)}
$$
Here $G(p,\rn)$ denotes the grassmannian of $p$-dimensional subspaces of $\rn$
and $\tr_WA$ denotes the trace of the quadratic form $A$ restricted to $W$.
Perhaps the simplest definition of $A\in \cp_p$ is given by requiring
$$
\l_1(A) +\cdots +\l_p(A)\ \geq0
\eqno{(\II.2)}
$$
where $\l_1(A) \leq \l_2(A) \leq \cdots \leq \l_n(A)$  are the ordered eigenvalues of $A$.

A function $u\in\USC(X)$ on an open set $X\ss\rn$ is said to be {\bf $p$-\psh on $X$}
if for each point $x\in X$ and each test function $\vf$ for $u$ at $x$, we have
$$
D_x^2 \vf \in \cp_p, \qquad
   {\rm i.e., }\ \ \tr_W\left\{ D^2_x \vf\right\} \ \geq\ 0\ \ \ \forall \, W\in G(p,\rn)
\eqno{(\II.3)}
$$
or equivalently
$$
\l_1(D_x^2 \vf)  +\cdots +\l_p(D_x^2 \vf)\ \geq0.
\eqno{(\II.3)'}
$$

The dual subequation $\wt{ \cp_p}$ can be defined by 
$$
\wt{ \cp_p}\ =\ \left\{  A\in \Symn : \tr_W \geq0 \ \ {\rm for\ some\ } W\in G(p,\rn)\right\}
$$
or equivalently, $A\in \wt{ \cp_p}$ if 
$$
\l_{n-p+1}(A) + \cdots + \l_n(A) \ \geq\ 0.
$$
Subharmonic functions for the dual subequation $\wt{ \cp_p}$ will be referred to as
{\bf dually $p$-plurisubharmonic functions}.  Function which are $\cp_p$-harmonic,
will be referred to simply as {\bf $p$-harmonic functions}.

The terminology {\sl $p$-\psh} is justified by the Restriction Theorem 7.3 in [\REST]
which says that $u\in\USC(X)$ is $p$-\psh if and only if the restriction of $u$ 
to each affine $p$-plane (or more generally to each connected $p$-dimensional  minimal submanifold) is subharmonic (or $\equiv -\infty$) for the classical
Laplacian (Laplace-Beltrami operator) on that $p$-plane.

The notion of $p$-plurisubharmonicity extends to any real number $p$ with $1\leq p\leq n$,
using the ordered eigenvalues as follows.
Let $\bar p = [p]$ denote the greatest integer in $p$.

\Def{\II.1}  For $A\in \Symn$ we say that 
$$
A\in\cp_p
\qquad\iff\qquad
\l_1(A) + \cdots + \l_{\bar p}(A) + (p-\bar p)\l_{\bar p+1} \ \geq\ 0
$$
%where $\l_1(A) \leq \cdots \leq \l_n(A)$ are the ordered eigenvalues of $A$.
\medskip
To see that $\cp_p$ is a convex cone, one shows that $\cp_p$ is
the polar of the set of symmetric forms $P_{e_1} +\cdots+ P_{e_{\bar p}} + (p-\bar p)
P_{e_{\bar p+1}}$ where $e_1,...,e_n$ are orthonormal.  Here $P_e$ denotes 
orthogonal projection onto the line through $e$.

If $\E$ is a subset of $\{u=-\infty\}$ where  $u\in \cp_p(X)$  (and $u$ is  not identically $-\infty$
on any component of $X$), then $\E$ is said to be {\bf $p$-pluripolar}.
Note that since $u\in \cp_p(X)$ implies that $\D u\geq0$,  the $p$-pluripolar sets can not have interior.
 Theorems \FF.1 and \GGG.1 lead to the following.

\Theorem{\II.2} 
{\sl
Suppose that $\E$ is a closed $p$-pluripolar set in $X$ and
$F\ss \Symn$ is any pure second order, constant coefficient subequation which is $\cp_p$-monotone.

\smallskip
(1)  Then $\E$ is removable for $F$-subharmonic functions which are locally bounded above across $\E$.
In particular, dually $p$-\psh functions on $X-\E$ extend to dually $p$-\psh functions on $X$ if they
are locally bounded above.\smallskip

(2) If $u\in C(X)$ is $F$-harmonic on $X-\E$, then $u$ 
 is $F$-harmonic on $X$. In particular, continuous functions on $X$ which are $p$-harmonic on $X-\E$ are 
 $p$-harmonic on $X$.
}

Since we have not assumed  $\psi$ to be smooth outside $\{\psi=-\infty\}$, 
the proof of Theorem \II.2 requires an extra step which is exactly like the one
given in the proof of Theorem \HH.1 above, and hence is omitted.

\Remark{\II.3}  Theorem \II.2 is vacuous if $1\leq p<2$ since there are no $p$-pluripolar sets $E$.
However, the classical ``$u+\e \psi$ technique'' employed in the proof can still be used with the
elementary Proposition \CC.3 replacing the fact that $u+\e \psi$ has no test functions at a 
point where $\psi$ equals $-\infty$.  See, for instance, the proof given in Appendix A of a
theorem of Caffarelli, Li and Nirenberg.

\medskip

Before describing lots of examples of subequations
$F$ which are $\cp_p$-monotone in Section \KK, we present some sufficient conditions
for a set $E$ to be $p$-pluripolar.

\vskip.3in

\noindent{\headfont \JJ.   $P$-Pluripolar Sets and Riesz Potentials.}
\medskip

The $p^{\rm th}$ {\bf Riesz kernel/function} $K_p$ is defined
(modulo a positive normalization) by 
$$
\begin{aligned}
K_p(x)\ =\   |x|^{2-p}, & \ \ 1\leq p<2,  \qquad K_2(x) \ =\ \log|x|, \ \ \ \  {\rm and}    \\ 
  &K_p(x)\  =\ -{1\over |x|^{p-2}}, \qquad    p>2.
\end{aligned}
\eqno{(\JJ.1)}
$$
They are fundamental for the subequation $\cp_p$.

\Prop{\JJ.2} {\sl
The Riesz kernel $K_p$ is $p$-harmonic on $\rn-\{0\}$ and $p$-\psh on $\rn$.
}

\pf
For $x\neq 0$ set $e=x/|x|$.  The functions $K_p$ have second derivatives $D^2_xK_p$ given (modulo a positive  constant multiple) by the single formula
$$
{1\over |x|^p}\left( P_{e^\perp} - (p-1) P_{e}\right)\ =\ 
{1\over |x|^p}\left(   I -  p  P_{e}\right)
\eqno{(\JJ.2)}
$$
This easily  implies that $K_p$ is $p$-harmonic on $\rn-\{0\}$.  Since  $K_p$ has no test functions at $x=0$, it is    $p$-plurisubharmonic on all of $\rn$.\qed

\medskip

Taking convex combinations of $K_p(x-y)$ via convolution yields a general class of $p$-\psh functions on $\rn$.

\Prop{\JJ.3} 
{\sl
Suppose $\mu$ is a compactly supported non-negative measure on $\rn$.  Then the 
$p^{\rm th}$ Riesz potential
$$
u\ \equiv\ K_p  * \mu
\eqno{(\JJ.3)}
$$
defines a $p$-\psh function on $\rn$ which vanishes at  $\infty$ for $p\geq3$ and is $\lloc(\rn)$.
}

\pf Replacing $K_p$ by $K_p^\a$, the maximum of $K_p$ and $-\a$, one obtains continuous functions
$u^\a \equiv K_p^\a * \mu$ which decrease down, as $\a\to\infty$, to the point-wise defined function
$u=K_p * \mu$.  Hence $u$ is upper semi-continuous. 
For integer $p$  its restriction to each affine $p$-plane $W$
is an (infinite) convex combination of $\D_W$-subharmonic functions, and as such is $\D_W$-subharmonic.
The proof for general $p$ is left to the reader.
\qed

\medskip

There is an extensive literature on Riesz potentials (see Landkof [\LAN]).
A compact set $\E\ss\rn$ is said to be {\bf $K_p$-polar} if it is contained in the 
$-\infty$ set of some Riesz potential $u=K_p * \mu$.   
In fact $\mu$ can be chosen so that $\supp \mu=\E$ and $\{K_p * \mu = -\infty\}=\E$
(see Section 1, Chapter III in [\LAN]).
Hence the Riesz potential is smooth in the complement of $\E$.
There is also a well defined notion of the 
$K_p$-capacity of $\E$ --  commonly called the $p-2$ capacity.  
Both of the facts:

\medskip
\centerline{
$E$ is $K_p$-polar $\quad\iff\quad$ $E$ has $K_p$-capacity (or $(p-2)$-capacity)  zero,
}

\medskip
\centerline{
$E$ has finite Hausdorff $(p-2)$-measure $ \quad\Rightarrow \quad $ 
the $K_p$-capacity of $E$ is zero.
}

%\medskip

\noindent
are classical (see Sections 1 and 4 of Chapter III in  [\LAN]).

Putting all this together with Proposition \JJ.3 we have the following.

\Prop{\JJ.4} 
{\sl
Suppose $\E$ is a compact subset of $\rn$.  Then:

\medskip
\centerline
{
 $\E$ has finite Hausdorff $p-2$-measure}
 
 \smallskip
 
 \centerline{
\quad$\Rightarrow$\quad $\E$ has $K_p$-capacity (i.e., $p-2$ capacity) zero}

\smallskip
\centerline{
\quad$\Rightarrow$\quad $\E$ is $C^\infty$ $p$-pluripolar. \qquad\qquad\qquad\qquad\qquad
}
}

\medskip

Further combining with Theorem \II.2 yields the following main result.

\medskip
\vfill\eject

\Theorem{\JJ.5} 
{\sl
Suppose $F\ss\Symn$ is a $\cp_p$-monotone subequation, and $E$ is a closed subset
of $X\ss \rn$.  If $E$ has locally finite Hausdorff $(p-2)$-measure, 
or, more generally, if $E$ has $p-2$ capacity zero, 
then $E$ is removable
for $F$-subharmonics and $F$-harmonics as described in parts (1) and (2) of Theorem \II.2.}

%%%%%%%%%%%%%%%%%%%%%%%%%%%%%%%%%%%%%%%%%%%%%%%
%%%%%%%%%%%%%%%%%%%%%%%%%%%%%%%%%%%%%%%%%%%%%%%
%%%%%%%%%%%%%%%%%%%%%%%%%%%%%%%%%%%%%%%%%%%%%%%
%%%%%%%%%%%%%%%%%%%%%%%%%%%%%%%%%%%%%%%%%%%%%%%
%%%%%%%%%%%%%%%%%%%%%%%%%%%%%%%%%%%%%%%%%%%%%%%

\vskip.3in
%\vfill\eject

\noindent{\headfont \KK.  Subequations which are   $\cp_p$-Monotone.}
\medskip

The previous two Theorems \II.2 and \JJ.5 require that the subequation $F$
be $\cp_p$-monotone.  Typically $F$ comes equipped with a natural monotonicity cone
$M$.  The purpose of this section is to characterize when $\cp_p\ss M$,
which ensures that $F$ is $\cp_p$-monotone.

\bigskip
\centerline{\bf Geometric Subequations.}
\medskip

First we illustrate with  a few examples,  starting with a large class -- the {\sl geometric} 
subequations $\M\equiv \cp(\GG)$  defined by fixing a closed subset $\GG\ss G(p,\rn))$ and setting
$$
\cp(\GG) \ \equiv\ \{A\in \Symn : \tr_W A \geq 0 \ \  \forall\, W\in \GG\}.
\eqno{(\KK.1)}
$$
Since $\GG \ss G(p,\rn)$, it is obvious that  $\cp_p \equiv \cp(G(p,\rn))$ is contained in $\cp(\GG)$.
Thus
$$
F \  {\rm is\ } \cp(\GG)\,{\rm monotone}
\qquad\Rightarrow\qquad
F \  {\rm is\ } \cp_p\,{\rm monotone}.
\eqno{(\KK.2)}
$$
This  means that 
$\cp(\GG)$ and its dual $\cpt(\GG)$ are $\cp_p$-monotone, and that Theorems \II.2 and \JJ.5 apply
to $\cpt(\GG)$.
These convex cone subequations have been studied in [\DDD] and  [\PUP], and in the more 
general riemannian setting in  [\DDR]. 
Here it is convenient to replace the term $\cp(\GG)$-{\sl subharmonic} 
by $\GG$-{\sl \psh}, and the term $\cp(\GG)$-{\sl polar}
by $\GG$-{\sl pluripolar}.

For an example, consider $\bbr^{2n}=\bbc^n$ and let $\GG\ss G(2,\bbr^{2n})$ be the set of complex lines.
Then $\cp(\GG)$ is exactly the subequation $\cp^\bbc$ discussed in Section \HH.

One also has $\GG =  {\rm LAG} \equiv \{W\in G(n,\bbr^{2n}) : W \ {\rm is\ Lagrangian}\}$.  There is a Lagrangian
equation of Monge-Amp\`ere type with many branches, each of which is $\cp(\LAG)$-monotone (see [\DDD]).
By Proposition \JJ.4 any set of locally finite Hausdorff (n-2)-measure in $\bbr^{2n}$  is $n$-pluripolar and therefore
LAG-pluripolar.

Another large class of examples comes from calibrations (see [\DDD]). Among the interesting
examples are the special Lagrangian, associative, coassociative and Cayley calibrations.
The Theorems  \II.2 and \JJ.5 apply  to each case.  For example, any set  $E\ss \bbr^8$ of locally 
finite Hausdorff 2-measure is Cayley pluripolar.  It would be interesting to find polar sets of larger 
dimension than that given by Proposition \JJ.4 for these   calibrations.  Preliminary 
work indicates that they may not exist. 

 \Remark {\KK.1}  It should be noted however that since $\cp(\GG)$ is a convex cone, the results of Section
 \LL \ apply to show that codimension-2 sets are removable for $\GG$-\psh
 functions which are locally bounded above (provided that $\GG$ involves all the variables).
Consequently, the new parts of Theorem \II.2 and Theorem \JJ.5 are part (1)  for $F=\wt\cp(\GG)$ (for
dually $\GG$-plurisubharmonic functions) and part (2) for $\GG$-harmonics.
\medskip

One can also consider quaternionic $n$-space $\bbr^{4n} = \bbh^n$ and
define $\cp^\bbh$ geometrically by the subset of quaternionic lines in $G(4,\bbr^{4n})$.
However,  $\cp^{\bbc, e} \ss \cp^\bbh$ for each of the complex structures defined on $\bbh^n$
by left multiplication by a unit imaginary quaternion $e$. 
In particular, this implies that the pluripolar sets for
each of these complex structures (e.g., the complex analytic  hypersurfaces) are
quaternionic polar.

\Ex{\KK.2.  (Branches of the $p$-Fold Sum Equation)}
For $p$ fixed, the convex cone subequation $\cp_p$ naturally occurs as a monotonicity subequation
for each of the following subequations.  In fact, $\cp_p$  is the first, or smallest, branch
of the family of subequations
$$
\cp_p \ =\  \cp_p^1\ \ss\   \cp_p^2\ \ss\  \cdots\ \ss\    \cp_p^N\ = \   \cpt_p 
$$
determined by the G\aa rding polynomial
$$
\MA_p(A)\ =\ \prod_{i_1<\cdots < i_p}\left( \l_{i_1} (A) +\cdots + \l_{i_p}(A)     \right)
$$
whose eigenvalues are the $p$-fold sums of the eigenvalues of $A\in \Symn$.
None of the other branches of $\MA_p$ are convex, but (see [\HYP] for the details)
$$
{\rm each\ branch\ of\ } \MA_p\ {\rm is\ } \cp_p\ {\rm monotone}.
\eqno{(\KK.3)}
$$

%The strongest kind of monotonicity for $F$ occurs by finding the largest monotonicity subequation.
%It is worth noting that the monotonicity subequations $\cp_k\equiv \cp(G(k,\rn))$ are nested
%$$
%\cp\ =\ \cp_1\ \ss\ \cp_2\ \ss\ \cdots\ \ss\ \cp_n,
%\eqno{(\KK.3)}
%$$
%as are  their complex and quaternionic analogues $\cp_k^\bbc\equiv \cp(G_\bbc(k,\rn))$,  i.e., 
%$$
%\cp^\bbc\ =\ \cp_1^\bbc\ \ss\ \cp_2^\bbc\ \ss\ \cdots\ \ss\ \cp_n^\bbc.
%\eqno{(\KK.4)}
%$$

\bigskip
\centerline{\bf The Riesz Characteristic.}
\medskip

The  following simple observation is the basis of this subsection.
Suppose  $F\ss \Symn$ is a subequation with monotonicity cone $M$.
Then if $\cp_p$ is contained in $M$, $\cp_p$ is also a monotonicity cone for $F$ and 
 Theorems \II.2 and \JJ.5 apply. 
 
 The cones $\cp_p$ are nested, that is, for real numbers $p<p'$ we have $\cp_p \ss\cp_{p'}$.
 A characterization of the largest $\cp_p$ contained in $M$ was given by 
 Theorem 5.1b in [\pPSH].  The statement, given below, uses the following easily computed
 invariant of $M$.

\Def{\KK.3}  The {\bf Riesz characteristic $p_M$ of $M$} is defined by
$$
p_M \ \equiv \ \sup \{ p : I-p P_e \in M \ \ {\rm for\ all\ } |e|=1\}.
$$
  
\Theorem{\KK.4. ($1<p<n$)}  {\sl
Suppose $M\ss\Symn$ is a convex cone subequation.  Then}
$$
\cp_p\ \ss\ M \qquad\iff\qquad p\ \leq\ p_M.
$$

\Cor{\KK.5} {\sl
The statement that $\cp_p\ss M$ is equivalent to either of the following:
\medskip

(1) \ \ $K_p(x)$ is $M$-subharmonic on $\rn$,

\medskip

(2) $(K_p* \mu)(x)$ is $M$-subharmonic on $\rn$ for all measures $\mu\geq0$ with compact support.
}

\pf Use (\JJ.2) and Proposition \JJ.3.\qed

\medskip

Because of Theorem \KK.4 the Theorems \II.2 and \JJ.5 can be recast in a more useful form
using the Riesz characteristic $p_M$ of $M$.

\Theorem{\KK.6}  {\sl
Suppose that $F\ss\Symn$ is a subequation with monotonicity cone $M\ss\Symn$
(a convex cone subequation). Let  $E\ss X$  be a closed subset.
Suppose any one of the following holds.

\medskip

(a) \ \ $E$ is $p_M$-polar.

\medskip

(b) \ \ $E$ has locally finite Hausdorff $(p_M-2)$-measure.

\medskip

(c) \ \ $E$ has $(p_M-2)$-capacity zero.
 
\medskip

\noindent
Then:
\medskip

(1)\ \ Each $u\in F(X-E)$ locally bounded above across $E$ has a canonical 
extension to $U\in F(X)$.
\medskip

(2)\ \ Each $u\in C(X)$ which is $F$-harmonic on $X-E$ is $F$-harmonic on $X$.

}

\bigskip

We conclude with several examples where this theorem applies. Specific cases are illustrated in Appendix B.

\Ex{1.  (The $\d$-Uniformly Elliptic Cone)}
For each $\d>0$ we define
$$
P(\d)\ \equiv\ \{ A\in \Symn : A+   \d ( \tr A) \cdot I \geq 0\}
%\eqno{(\KK.6)}
$$
Note that 
$$
I-p P_e \in \cp(\d) 
\qquad\iff\qquad
I-p P_e +\d(n-p)\geq0
\qquad\iff\qquad
p\ \leq\ {1+\d n\over 1+\d}.
$$
This proves that
$$
\cp(\d) \ \ {\rm has\ Riesz\  characteristic}\ \  p\ =\ {1+\d n\over 1+\d}.
$$
It is shown in [\ASPECTS, Lemma A.1] that all other O$(n)$-invariant convex cone subequations with 
Riesz characteristic $p$ are contained in this $\cp(\d)$.

\Ex{2.  (The Pucci Cone)} For $0<\l <\L$ we define
$$
\cp_{\l,\L} \ \equiv \ \{A\in \Symn : \l\tr A^+ + \L\tr A^- \geq0\},
$$
where $A=A^+ + A^-$ is the decomposition into $A^+\geq0$ and $A^-\leq0$.
Since $I-pP_e = P_{e^\perp} - (p-1)P_e$, we have that
$$
\cp_{\l,\L} \ \ {\rm has\ Riesz\  characteristic}\ \ {\l\over \L}(n-1) +1.
$$

The condition that a subequation $F$ be uniformly elliptic (in the standard sense)
can be restated as ``monotonicity'' using either of the two cones $\cp(\d)$ or $\cp_{\l,\L}$.
That is, $F$ is uniformly elliptic if and only if either
$$
\begin{aligned}
&F+\cp(\d) \ \ss\ F \ \ \ {\rm for\ some}\ \ \d>0, {\rm or}  \cr
&F+\cp_{\l, \L} \ \ss\ F \ \ \ {\rm for\ some}\ \  0<\l<\L.\cr
\end{aligned}
$$
(See Section 4.5 in [\SURVEY] for more details.)

\Ex{3.  (The $k^{\rm th}$ Elementary Symmetric Cone)}
This equation is often referred to as the $k^{\rm th}$ hessian equation.
 For $k=1,...,n$ we define
$$
\Sigma_k \ \equiv \  \{A \in \Symn : \s_1(A)\geq0, ... , \s_k(A)\geq0 \},
$$
where $\s_j(A)$ is the $j^{\rm th}$ elementary symmetric function in the eigenvalues of 
$A$.  One easily computes that 
$$
\Sigma_k \ \  {\rm has\ Riesz\  characteristic}\ \ { n \over k}.
$$

In all three of the examples above there is a  polynomial equation on $\Symn$ which is 
G\aa rding hyperbolic with respect to the identity (see [\HYP]),
and the cone $M$ is the G\aa rding cone associated with this  polynomial.
The G\aa rding polynomial determines other natural subequations called
branches which are $M$-monotone.

Theorem \KK.6 can be applied to a broad spectrum of nonlinear equations.
Here are some cases related to the examples above.

\Ex{1$'$} \ If $E$ has locally finite Hausdorff ${\d(n-2)-1 \over \d+1}$ measure, 
then $E$ is removable (as in the conclusion of Theorem \KK.6) for each branch
of the equation $\det(A+\d(\tr A) I) =0$.
\smallskip

In fact, this  removable singularity result  holds for the 
{\bf $\d$-elliptic regularization} 
$
F(\d) \ \equiv \ \{A\in\Symn : A+\d(\tr A)I \in F\}
$
 of {\sl any} pure second-order subequation $F$, because $F(\d)$ is always 
 $\cp(\d)$-monotone.
This follows from the general fact that 

\smallskip
\centerline
{\sl if $M$ is a monotonicity cone for $F$,
then $M(\d)$ is a monotonicity cone for $F(\d)$.
}

\Ex{2$'$} \ If $E$ has locally finite Hausdorff ${\l\over \L}(n-1)-1$ measure, then $E$
is removable (as in the conclusion of Theorem \KK.6) for each branch
of the G\aa rding equation  defining $\cp_{\l,\L}$.
\smallskip

This G\aa rding polynomial is defined as follows. For each subset $I\ss\{1,...,n\}$
define $v(I)\in \rn$ by $v(I)_i =\l$ if $i\in I$ and $v(I)_i =\L$ if $i\notin I$. The points   
$v(I)$ are the vertices  of the cube $[\l, \L]^n\ss\rn$.  Let ${\cal I}$ denote the subset
of $I\ss\{1,...,n\}$ such that the open segment from 0 to $v(I)$ is disjoint from the cube.
Given $A\in\Symn$, let $\l_I(A) \equiv \sum_{i\in I}\l_i(A)$. Finally, set $I' \equiv
\{1,...,n\} - I$. The the G\aa rding polynomial with G\aa rding cone $\cp_{\l,\L}$ is
$$
p(A) \ \equiv\ \prod_{I\in {\cal I}} \left( \l \, \l_I(A) + \L \, \l_{I'}(A)\right).
$$
Note that with $\L = 1+\d$ and $\l = \d$, the degree $n$ polynomial $\det(A+\d(\tr \,A)I)$
defining $\cp(\d)$ (as its G\aa rding cone) is a factor of this polynomial $p(A)$ defining $\cp_{\l,\L}$.

\Ex{3$'$} \   If $E$ has locally finite Hausdorff ${n \over k}-2$ measure, then $E$
is removable (as in the conclusion of Theorem \KK.6) for each branch
of the equation  $\s_k(A)=0$.

\bigskip
\centerline{\bf Variable Coefficients.}
\medskip

The ideas above can be applied directly to  variable coefficient subequations
$F\ss \bbj^2(X) \equiv X\times \bbr\times \rn\times \Symn$ as defined in [\DDR]. Suppose that such 
a subequation $F$ on an open set $\O\ss\rn$ is ${\bf P}_p$-monotone
where ${\bf P}_p \equiv \bbr\times\rn \times \cp_p$ (i.e., $F+{\bf P}_p\ss F$ under  fibre-wise sum). 

 \medskip
\centerline
{\sl
Then any closed set $E\ss\O$ of locally finite Hausdorff $(p-2)$-measure }
\centerline
{\sl
 is removable for $F$-subharmonics and $F$-harmonics as in Theorems 6.1 and 7.1.}

\pf Consider a Riesz potential $\psi$ which is $-\infty$ on $E$ and
 $\cp_p$(i.e., ${\bf P}_p$)-subharmonic as above.
%Then $\psi$ is, of course,  ${\bf P}_p$-subharmonic, and so 
Then $E$ is $C^\infty$ ${\bf P}_p$-polar.
The result then follows from Theorems 6.1 and 7.1.
\qed

\medskip

Now on any riemannian manifold there is a natural subequation ${\bf P}_p^X$ which
generalizes ${\bf P}_p$ on euclidean space. It is defined as in Definition 9.1 by using
the ordered eigenvalues of the riemannian hessian (cf. \S 4 in [\DDR]).

\Theorem{\KK.7} {\sl
Let $F$ be a subequation on a riemannian manifold $X$ which is ${\bf P}_p^X$-monotone.
Then any closed subset $E\ss X$ of Hausdorff dimension $p'<p$  
is removable for $F$-subharmonics and $F$-harmonics as Theorems 6.1  and 7.1
}

\pf
Fix $x\in X$ and choose geodesic normal coordinates at $x$.  Straightforward calculation
of $D^2K_p$ shows that for  any fixed $p'<p$ there exists $\e>0$ so that 
$\Hess_x\{K_{p'}(x-y)\} \in \cp_p$ (i.e. $K_{p'}(x-y)$ is ${\bf P}_p$-subharmonic) 
for all $x,y \in B(0,\e)$. Thus if $\mu$ is a measure with support in $B(0,\e)$, then
$\mu * K_{p'}$ is ${\bf P}_p$-subharmonic on $B(0,\e)$. Therefore, assuming that 
$E \cap B(0,\e)$ has Hausdorff dimension $\leq p'$ we conclude that $E \cap B(0,\e)$
is $C^\infty$ ${\bf P}_p$-polar, and Theorems 6.1 and 7.1 apply.\qed

%Then in these coordinates $({\bf P}_p^X)_0=
%({\bf P}_p)_0$.  It follows that for any $p'<p$ there exists $\e>0$ so that the constant
%coefficient subequation ${\bf P}_{p'} \ss {\bf P}_p^X$ on the $\e$-ball $B(0,\e)$.
%The theorem now follows from the result above.\qed

%%%%%%%%%%%%%%%%%%%%%%%%%%%%%%%%%%%%%%%%%%%%%%%%%
%%%%%%%%%%%%%%%%%%%%%%%%%%%%%%%%%%%%%%%%%%%%%%%%%
%%%%%%%%%%%%%%%%%%%%%%%%%%%%%%%%%%%%%%%%%%%%%%%%%
%%%%%%%%%%%%%%%%%%%%%%%%%%%%%%%%%%%%%%%%%%%%%%%%%
%%%%%%%%%%%%%%%%%%%%%%%%%%%%%%%%%%%%%%%%%%%%%%%%%
%%%%%%%%%%%%%%%%%%%%%%%%%%%%%%%%%%%%%%%%%%%%%%%%%

\vskip.3in
%\vfill\eject

\noindent{\headfont \LL.   Removable Singularities for Convex Subequations.}
\medskip

The results of this section do not depend on Theorems \FF.1 and \GGG.1.

A Theorem stated in [\HAR, p.132], and again in [\HPP, Thm. 1.2], for classical \psh functons $u$ 
(the case where $F=\cp^\bbc$ as in Section \HH)
contains four removable singularity results which entail different local growth hypotheses on $u$ across
$E$, namely:

\smallskip
\centerline
{
\qquad
(a)\ \ $u\in C^\a(X)$, \ ($0<\a\leq 1$), \qquad \qquad \qquad\ \  (b)\ \ $u\in L^p_{\rm loc}(X)$, \ ($1\leq p < \infty$) \qquad\ \ 
}

\smallskip
\centerline
{
 (c)\ \ $u$ locally bounded above across $E$, \ \ \qquad(d)\ \ no condition on $u$ across $E$.
}

\smallskip
\noindent
The different size hypotheses on $E$ are:

\smallskip

(a) \ \ $E$ has zero Hausdorff $2n-2+\a$ measure.

\smallskip

(b) \ \ $E$ has locally finite Hausdorff $2n-2q$ measure, with ${1\over p}+{1\over q}=1$.

\smallskip

(c) \ \ $E$ has locally finite Hausdorff $2n-2$ measure.

\smallskip

(d) \ \ $E$ has zero  Hausdorff $2n-2$ measure.

\smallskip

Parts (a) and (b) are due to the first author and Polking [\HP], [\HPP], part (c) is due to Lelong [\LEL],
and part (d) is due to Shiffman [\SHF].

Each of the first three results in this theorem has a straightforward generalization 
to subharmonics (but not harmonics) for any  convex subequaton $F$ which is second-order complete.
A subequation $F$ is {\sl convex} if each fibre $F_x$ is convex.
It is  {\sl  second-order complete} if each non-empty fibre $F_{x,r,p} \ss\Symn$
``depends on all the variables in $\rn$''.  This means that there does not exist a proper subspace
$W\ss\rn$ and a subset $F'\ss\Sym(W)$ such that
$$
A\in F_{x,r,p}\quad\iff\quad A\bigr|_W \in F'.
$$
(When $F_x$ is convex, if one non-empty fibre $F_{x,r,p} \ss\Symn$
depends on all the variables in $\rn$, then so does every other non-empty
fibre  $F_{x,r',p'}$.)

%(Here $ A\bigr|_W$ denotes the restriction of the quadratic  form to $W$.)

Throughout this section we make both of  the above assumptions on $F$ -- convexity and 
  second-order completeness --   in addition to the 
positivity condition (P) and  the negativity condition (N) (see Definition 3.8 in [\DDR]).

\Theorem { \LL.1}
{\sl
Suppose $E$ is a closed subset of an open set $X\ss\rn$, and $u$ is $F$-subharmonic on $X-E$.

\medskip
(a) \ \ If $u\in C^\a(X)$ and  $E$ has  Hausdorff $(n-2+\a)$-measure   zero,
 then $u$ is $F$-subharmonic on $X$. \medskip

(b) \ \ If $u\in {\rm L}_{\rm loc}^p(X)$  ($1\leq p<\infty$)\   and  $E$ has  locally finite Hausdorff  $(n-2q)$-measure  
(${1\over q}+{1\over p}=1$), 
 then $u$ is $F$-subharmonic on $X$. 

 \medskip
(c) \ \ If $u$ is locally bounded above across  $E$ and  $E$ has Newtonian capacity zero,
 then the canonical extension $U$ of $u$ to $X$ is $F$-subharmonic.
 In particular,  if  $E$ has  locally finite Hausdorff  $(n-2)$-measure, then the capacity of $E$ is zero
and hence $U\in F(X)$.
 }
 \bigskip
 
The fourth part (d)  of the theorem in [\HAR], namely:

 \medskip
 \centerline{ If $u$ is \psh on $X-E$ and $E$ has Hausdorff $(2n-2)$-measure zero,}
 
 \centerline{
 then $u$ has a \psh extension to $X$. }
\medskip

\noindent
 has the following counterpart in $p$-geometry (see Sections \II \ and \JJ\  above.)

\Theorem{\LL.2. ($p\geq2$)} 
{\sl
Suppose $u$ is $p$-\psh on $X-E$ and the Hausdorff $(n-p)$-measure of the closed set $E$ is zero.
Then $u$ has a $p$-\psh extension to $X$
}

\pf
This is an immediate consequence of Part (c) of Theorem  \LL.1 and  Lemma \FF.2.\qed

\medskip
\noindent
{\bf Discussion of the Proof of Theorem  \LL.1.} The proof is by reduction to the linear case.
\medskip

\medskip
\noindent
{\bf Step 1. (The Strong Bellman Principle).}
The previous assumptions on $F$ (together with a mild regularity assumption)
 imply that locally there exists a family $\cf$ of linear (sub)equations
of the form
$$
Lu\ =\ \bra {a(x)}{D_x^2 u} + \bra {b(x)}{D_x u}  - c(x) u \ \geq\ \lambda
\eqno{( \LL.1)}
$$
with $a(x)>0$ positive definite and $c(x)\geq0$ at each point $x$
with the property that 
$$
\begin{aligned}
&{\rm A \ function \ } u \ {\rm  is \ }\, F\,{\rm subharmonic} \quad \iff \quad \cr
{\rm locally\ } &u {\rm\  satsfies\ } Lu \geq \l\ {\rm  for\  all\  pairs\ }
(L,\l)\in \cf.
\end{aligned}
\eqno{( \LL.2)}
$$
This is Theorem 8.6 in  [\BELL].
\medskip

\medskip
\noindent
{\bf Step 2.}
Show that a viscosity $(L,\l)$-subharmonic function $u$ belongs to $\lloc$ and
satisfies $Lu\geq \l$ in the distributional sense.

\medskip
\noindent
{\bf Step 3.}
Prove the theorem in the linear distributional case.  For solutions to $Lu=\l$,
Part (a) is Theorem 4.4 in [\HP] while 
part (b) is Theorem 4.1(a) in [\HP].  For proving part (a) the modifications required to treat
$Lu\geq\l$ are described on pages 705-706 in [\HPP], while for part (b) the modifications
for $Lu\geq\l$, as described on pages 132-133 in [\HAR], go as follows. 

Pick  $\psi$ smooth and compactly supported in $X$ with $\psi\geq0$.
For convenience assume $\int\psi = 1$.  One must show that $(Lu)(\psi)\geq \int \l \psi$ for all such 
$\psi$.  Choose $\vf_\e$ as in Lemma 3.2 in [\HP]. Now $(Lu)(\psi) = 
(Lu)((1-\vf_\e)\psi)   +(Lu)(\vf_\e\psi)$. The estimate on the derivatives of $\vf_\e$
provided by Lemma 3.2 in [\HP] show that $L^t(\vf_\e\psi)$ converges to zero in $L^q$, 
and hence by H\"older's inequality that $(Lu)(\vf_\e\psi)$ converges to zero.
Meanwhile, $(Lu)((1-\vf_\e)\psi)  \geq  \int \l (1-\vf_\e)\psi$ (since $Lu\geq\l$ on $X-E$)
implies that $(Lu)(\psi) \geq \int\l\psi$.

Perhaps it is worth noting that neither of  the conditions $L$ elliptic ($a(x)>0$) or $c(x)\geq0$ 
was used to conclude that $Lu\geq\l$ across $E$ in the distributional sense.  In fact, the only thing
 used about  $L$ was that it has order 2.

Part (c) in the linear distributional case is classical.
If $E$ has Newtonian capacity zero, then $E$ is $L$-polar.  Therefore our
Theorem \FF.1 applies with $F$ defined by $Lu\geq\l$ and $M=L$.  Of course our proof is just the
classical proof of the removability of $E$ for $F$-subharmonic functions.

\def \esssup#1{{\rm ess} \!\! \sup_{#1}}
\def \essup#1{{\rm ess}  \sup_{#1}}

\medskip
\noindent
{\bf Step 4.}
Show that distributional solutions $u$ to $Lu\geq \l$ (they belong to $\lloc$) have a pointwise
canonical representative $U\in\USC(X)$ in the $\lloc$-class given by
$$
U(x)\ =\  {\rm ess \, \overline{ \lim_{y \to x}}} \ u(y) \ =\   \lim_{\rho\to0} \esssup {B_\rho(y)} \ u
$$
The important point here is that $U$ is independent of the operator $L$,
and that $U$ is $(L,\l)$-subharmonic in the viscosity sense.
(See the Appendix in [\AC] and Section 9 in [\BELL].)

\medskip
\noindent
{\bf Step 5.}
Use (\LL.2), this time in the reverse direction, to conclude that $U$ is $F$-subharmonic on $X$.

%%%%%%%%%%%%%%%%%%%%%%%%%%%%%%%%%%%%%%%%%%%
%%%%%%%%%%%%%%%%%%%%%%%%%%%%%%%%%%%%%%%%%%%
%%%%%%%%%%%%%%%%%%%%%%%%%%%%%%%%%%%%%%%%%%%
%%%%%%%%%%%%%%%%%%%%%%%%%%%%%%%%%%%%%%%%%%%

%\vfill\eject
\vskip .3in

\centerline{\headfont  Appendix A.  }\medskip

\centerline{\headfont  A Removable Singularity  Theorem}

\centerline{\headfont  of Caffarelli, Li and Nirenberg}
\medskip

The classical ``$u+\e\psi$-technique'', that we have used to prove Theorem \FF.1, 
naturally lends itself to prove   a recent result of 
Caffarelli, Li and Nirenberg [\CLN, Thm. 1.3].    
Let $d_\E(x) \equiv \dist(x,\E)$ denote the distance to a subset  $\E\ss\rn$.

\Theorem {A.1. (Caffarelli, Li and Nirenberg)} 
{\sl
Suppose $\E$ is a closed submanifold of an open set $X$ in $\rn$, and that $F$
is a general subequation on $X$.  Consider a function $u\in F(X-\E)$ 
which is locally bounded above across $\E$.   Let $U\in \USC(X)$ be the
canonical upper semi-continuous extension of $u$ to $X$.  If for all $\e>0$ sufficiently small
$$
U+\e d_\E\ \ {\sl has\ no\  test\ functions\ at\ points\ of\   }\E,
\eqno{(A.1)}
$$
then $U$ is $F$-subharmonic on $X$.
}

\medskip

In [\CLN] a function $U\in\USC(X)$ satisfying (A.1) is said to be {\bf upper conical on $\E$}.

\pf
Note that $U+\e d_\E$ decreases to $U$ as $\e\searrow 0$.  Therefore, by the 
{\sl Decreasing Limit Property" }(see, for example, Theorem 2.6 in [\DDR]), if each
$U+\e  d_\E$ is $G$-subharmonic on $X$  for some fixed subequation $G$, then $U$ is also
$G$-subharmonic on $X$.  Now for each $c>0$ take $G$ to be  the enlargement $F^c$ of our subequation 
$F$ defined by 
$$
F^c_x \ \equiv\ \{J = (r,p,A) : \dist_x(J,F_x)\leq c\},
\eqno{(A.2)}
$$
where the fibre-distance  is defined so that $F^c = F+S^c$ for the constant 
coefficient subequation
$
S^c \equiv (-\infty, c] \times B_c(0)\times (\cp-c\cdot {\rm I}).
$
It will suffice to show that: for each $c>0$
$$
U+\e d_\E \in F^c(X)\ \ {\rm for\ all\  } \e\ {\rm sufficiently\ small},
\eqno{(A.3)}
$$
since the assertion: $U\in F^c(X) \ \forall \, c>0 \ \Rightarrow\ U\in F(X)$  is a trivial
consequence of the definition  of subharmonicity.   This uses the fact  that $F= \bigcap_{c>0} F^c$
which follows from  positivity (P) and negativity (N) for $F$.

By the hypothesis (A.1) and the trivial Proposition \CC.3 it will suffice to prove 
$$
U+\e d_\E \in F^c(X-E)\ \ {\rm for\ all\  } \e\ {\rm sufficiently\ small}.
\eqno{(A.3)'}
$$
In fact, since the result is local, we may replace the submanifold $\E$ by the interior of a 
small compact subset of $\E$ and replace $X$ by a thin normal neighborhood
 of $\E$.

Given a point $x\notin \E$ but near $\E$, there exists a unique line segment from a point
$x_0\in \E$ to $x$, of length $\d=d_E(x)$. We write $x = x_0+\d \nu$. 
Let $B$ denote the normal disk to $\E$ of radius $\d$ at $x_0$.  The tangent space   at $x$ splits as
$$
T_x\rn \ =\ T_x(B)^\perp\oplus T_x(\partial B) \oplus \span \nu.
$$
Let $II_x^\E$ denote the second fundamental form of $\E$ in the normal direction $\nu$, 
translated along the line segment from
$x_0$ to $x$ and acting as a quadratic form on $T_x(B)^\perp$.  
Let $P_W\in\Symn$
denote orthogonal projection onto a subspace $W$.  By direct computation the 
2-jet of $d_\E$ is given by:
$$
d_\E(x) = \d,\ \ D_xd_\E = \nu,\ \ D^2_x d_\E = II_x^\E + {1\over \d} P_{T_x(\partial B)}.
\eqno{(A.4)}
$$

We can assume that $d_E(x)<1$ for all $x\in X$.
Note that $|D_x d_E|=1$.
Finally, $D^2_x d_\E\geq -\kappa I$
 where $\kappa = \sup |\kappa(\nu')|$  
taken over all principal curvatures of all  normal directions to $\E$ (which is finite
by our pre-compactness assumption). Thus with $C\equiv \max\{1,\kappa\}$ we have
$$
\left(d_E(x), D_x d_E, D^2_x d_E    \right) \ \in \ S^C \ \ \forall\, x \in X-E,
\eqno{(A.5)}
$$
that is, $d_E$ is $S^C$-subharmonic on $X-E$.
It follows directly from Definition \BB.1,  first, that $\e d_\E$ is 
$S^{\e C}$-subharmonic on  $X-E$, and second, since $d_E$ is smooth, that $U+ \e d_\E$
is 
$F^{\e C}$-subharmonic on $X-E$.
This proves (A.3)$'$ if $\e<1$ since
 $F^{\e C}(X-E)\ss F^{C}(X-E)$ for $0<\e<1$.
\qed

\medskip
\def\xb{\bar x}

Remark 1.2 in [\CLN] can be used to replace the hypothesis (A.1) on $U+\e d_E$
with a hypothesis on $U$ itself.
The only way the hypothesis (A.1) on $U+\e d_E$ can fail to be true at $\bar x\in E$ is for $U$
to have a very large 2-jet at $\bar x$.  More precisely: 

\medskip
\centerline
{\sl  if $U+\e d_E$ has a test function 
$\vf$ at $\bar x$ and if $\psi$ is any test function for $-\e d_E$ at $\xb$,}

\centerline{\sl  then  $\vf+\psi$ is a test function for $U$ at $\xb$.}  

\medskip

This can be restated in terms of jets  as follows.
For an upper semi-continuous function $v$ defined near $\xb$ we let
$$
J_{\xb}^{+} v   \ = \  \{(D_{\xb} \vf, D^2_{\xb}) : \vf \ {\rm is\ a\ test\ function\ for\ } v \ {\rm at\ } \xb\}
$$
denote the {\sl upper (reduced) 2-jet}  of $v$ at $\xb$.  
$$
(A.1) \ {\rm is\ false \ at \ } \xb, {\rm \ i.e.,\ } J_{\xb}^{+} (U +\e d_E) \neq \emptyset \qquad \Rightarrow
$$
$$
J_{\xb}^{+}  (U +\e d_E)  + J_{\xb}^{+}  (-\e d_E) \ \ss\ J_{\xb}^{+}  (U).
\eqno{(L)}
$$
Thus any condition which limits the size of the 2-jet of $U$ at $\xb$ enough to violate the 
``large'' condition (L) will imply the hypothesis (A.1).
For example when $E$ is a point, it is an easy calculation to see that
$$
J_{\xb}^{+} (-\e|x-\xb|) \ =\ B_\e(0)\times \Symn.
\eqno{(A.6)}
$$
Thus if (A.1) fails when $E=\{\xb\}$, then
$$
B_\e(p) \times \Symn\ \ss\  J_{\xb}^{+}  (U).
\eqno{(L')}
$$
where $p$ is the first derivative of a test function for $U+\e|x-\xb|$  at $\xb$.

 \Cor{A.2}  
 {\sl
 Suppose $u\in F(X-\{\xb\})$ with $u$ bounded above across $\xb$.  Then the extension
 $U \in\USC(X)$ of $u$, given by defining $U(\xb)=\overline{\lim}_{x\to\xb} u(x)$, is $F$-subharmonic
 on $X$, unless the reduced upper 2-jet of $U$ at $\xb$ contains $B_\e(p)\times \Symn$
for some $p\in \rn$ and $\e>0$.

 }

%%%%%%%%%%%%%%%%%%%%%%%%%%%%%%%%%%%%%%%%%%%
%%%%%%%%%%%%%%%%%%%%%%%%%%%%%%%%%%%%%%%%%%%
%%%%%%%%%%%%%%%%%%%%%%%%%%%%%%%%%%%%%%%%%%%
%%%%%%%%%%%%%%%%%%%%%%%%%%%%%%%%%%%%%%%%%%%

%\vfill\eject
\vskip .3in

\centerline{\headfont  Appendix B. An Illustration:}
\medskip
\centerline{\headfont  Cones with $p={3\over2}$ and $n=3$  }\medskip

\medskip

\centerline{
{\includegraphics[width=.36\textwidth]{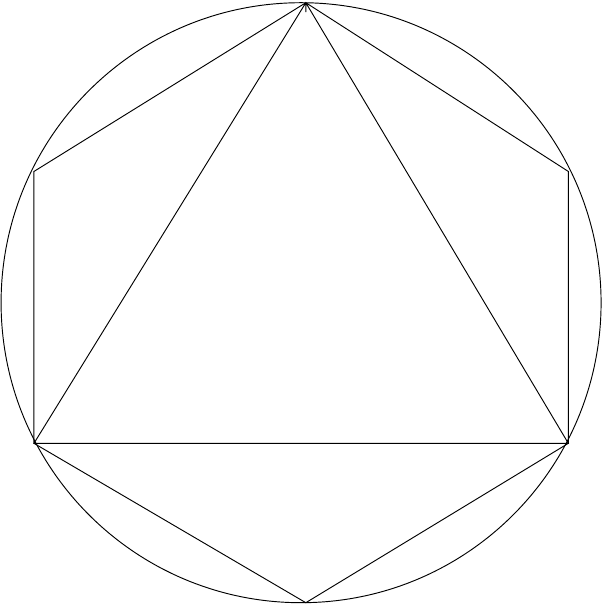}} \qquad
{\includegraphics[width=.4\textwidth]{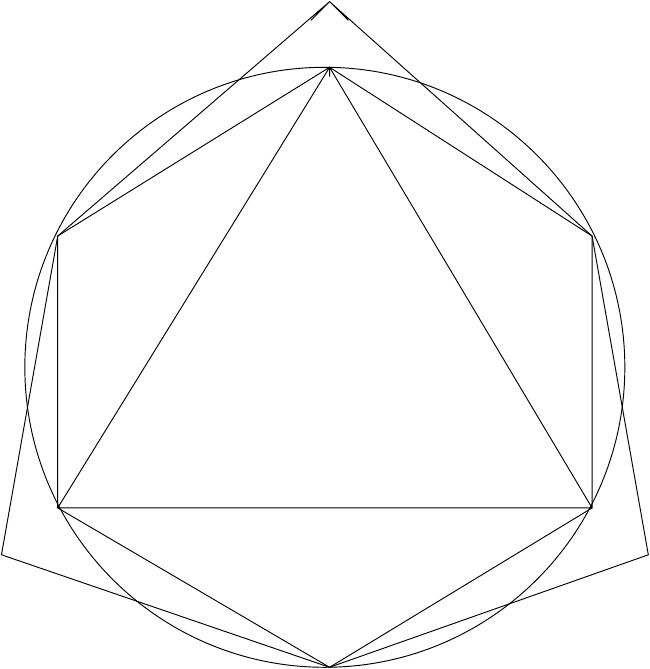}}
}

\centerline
{
$\cp\ \ss\ \cp_{3\over 2} \ \ss\ \Sigma_2$
\hskip1in
Pucci Cone Added
}

\centerline{
{\includegraphics[width=.7\textwidth]{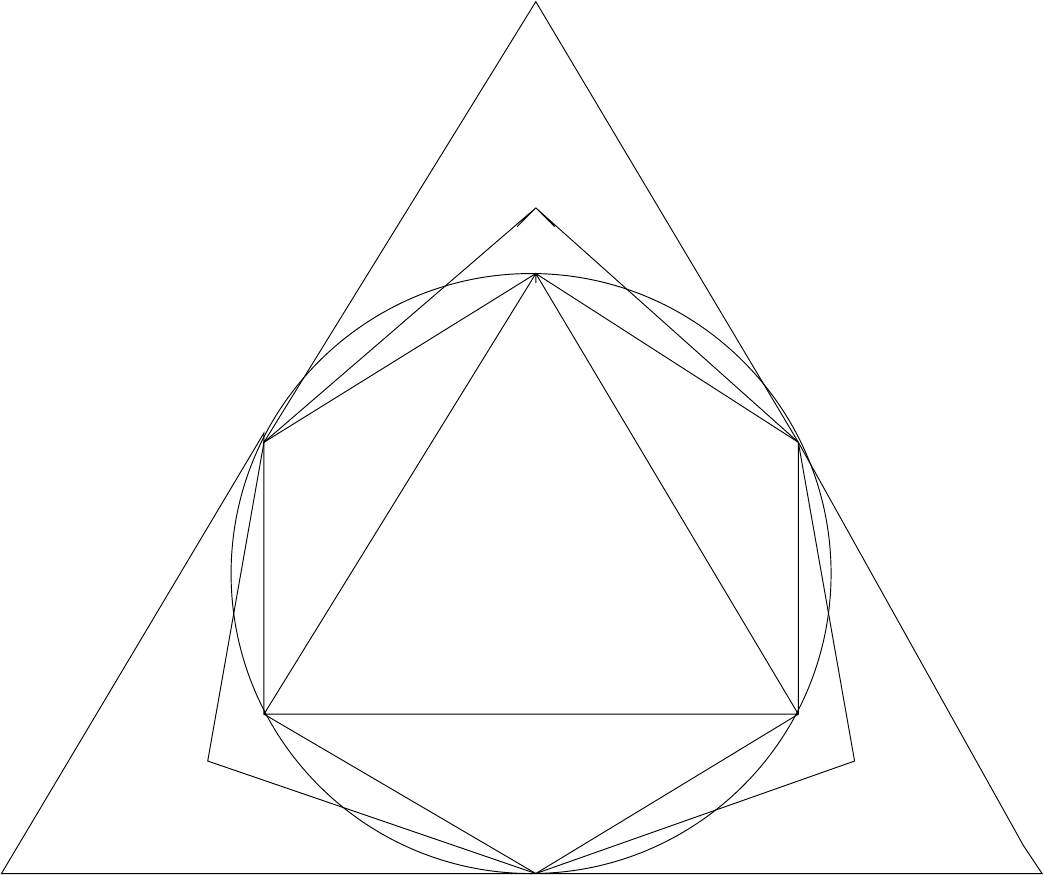}}
}

\smallskip
\centerline
{
$\cp(\d)$ Added
}

These  figures illustrate various convex cones discussed in the
paper in the specific case of $\Sym(\bbr^3)$ and Riesz characteristic
$p={3\over 2}$.  Since  the sets are
SO(3)-invariant, it suffices to represent them in eigenvalue space
$$
\bbr^3\   =\ \{(\l_1,\l_2\l_3) : \l_k\in\bbr\}. 
$$
Moreover, it suffices to
consider the bases of these cones in
$$
\bbr^2\   =\ \{(\l_1,\l_2\l_3) \in \bbr^3 : \l_1+\l_2+\l_3=1\}.
$$
 
There are two reference cones:  the cone $\cp$ whose base is the inscribed triangle,
and the cone $\D\equiv\{\tr A>0\}$ whose base is entire plane.  All other cones represented
have Riesz characteristic ${3\over2}$. In monotone progression they are as follows.

\medskip

{\bf The Inscribed Hexagon} is the base of the  cone $\cp_{3\over2}$.   For slightly larger values of $p$, the base of the cone $\cp_p$ is also an (irregular) hexagon.  In these cases, three of the vertices coincide with the vertices of the triangle.  The other
three vertices are further from the origin than those of $\cp_{3/2}$.

\medskip

{\bf The Disk}  (which has radius $\sqrt{{3\over2}}$)  is the base of the  second Hessian cone $\Sigma_2$, given by  the conditions $\s_1(A)\geq0$ and $\s_2(A)\geq0$. 
(See Example 3 in Section \KK.)
By the remarks  above, $p={3\over 2}$
is the largest $p$ with $\cp_p\ss\Sigma_2$.

\medskip

{\bf  The Second Hexagon} is the base of the Pucci  cone $\cp_{\l,\L}$  of Riesz characteristic
${3\over2}$.   (See Example 2 in Section \KK.) It is the smallest Pucci which contains 
$\cp_{3\over2}$. This holds if and only if ${\l\over\L}= {1\over4}$. Note that the base of $\cp$ is contained
in the {\sl interior} of the Pucci cone.  This is equivalent to uniform ellipticity of the Pucci equations.

\medskip

{\bf The Circumscribed Triangle}  is the base of the uniform ellipticity   cone $\cp(\d)$   of Riesz characteristic ${3\over2}$.  (See Example 1 in Section \KK.) Here $\d={1\over 3}$. The bases of the
family of cones $\cp(\d)$ differ by homotheties from the center of the disk. The circumscribed 
triangle is clearly the smallest such to contain $\cp_{3\over2}$.

 \medskip
 
 Note the three points common to the boundaries of each of the four subequations with common characteristic $p={3\over 2}$.
 They represent ${2\over3}I-P_{e_j}$, $j=1,2,3$ (cf. Definition \KK.3).

\def\item{}
%\vfill\eject
\vskip .3in

\centerline{\bf References}

\vskip .2in

 \smallskip

\noindent
\item{[\ALEX]}     A. D. Alexandrov,  {\sl  The Dirichlet problem for the equation Det$\| z_{i,j}\| = \psi(z_1,...,z_n,x_1,...,x_n)$}, I. Vestnik, Leningrad Univ. {\bf 13} No. 1, (1958), 5-24.

\noindent
\item{[\AGV]} M. E. Amendola, G. Galise and A. Vitolo, 
{\sl Riesz capacity, maximum principle and removable sets of fully
nonlinear second-order elliptic operators},  Preprint, University of Salerno.

\smallskip

 \noindent
\item{[\ASSA]}  S. N. Armstrong,  B. Sirakov and C. K.  Smart,  {\sl Fundamental
solutions of homogeneous fully nonlinear elliptic equations}, Comm.
Pure. Appl. Math. {\bf 64}  (2011), no. 6, 737-777.

\smallskip

 \noindent
\item{[\ASSB]}  \ \----------,  {\sl Singular solutions
of fully nonlinear elliptic equations and applications}, Arch. Ration.
Mech. Anal. {\bf 205} (2012), no. 2, 345-394.

\smallskip

\noindent
\item{[\BTA]}   E. Bedford and B. A. Taylor,  {The Dirichlet problem for a complex Monge-Amp\`ere equation}, 
Inventiones Math.{\bf 37} (1976), no.1, 1-44.

\smallskip

\noindent
\item{[\BTB]}  \ \----------,   {Variational properties of the  complex Monge-Amp\`ere equation, I. Dirichlet principle},     Duke  Math. J. {\bf 45} (1978), no. 2, 375-403.

\smallskip

\noindent
\item{[\BTC]}   \ \----------,   {A new capacity for plurisubharmonic functions}, 
Acta Math.{\bf 149} (1982), no.1-2, 1-40.

\smallskip

 \noindent
\item{[\CLN]}  L.  Caffarelli,  Y.Y.  Li and L. Nirenberg {\sl  Some remarks on singular solutions
of nonlinear elliptic equations, III: viscosity solutions including parabolic operators},
  Comm. on Pure and Applied Math. {\bf 66} (20013),  109-143.  ArXiv:1101.2833.

\smallskip

\noindent
\item{[\CIL]}   M. G. Crandall, H. Ishii and P. L. Lions {\sl
User's guide to viscosity solutions of second order partial differential equations},  
Bull. Amer. Math. Soc. (N. S.) {\bf 27} (1992), 1-67.

 \smallskip

\noindent
\item{[\CRA]}   M. G. Crandall,  {\sl  Viscosity solutions: a primer},  
pp. 1-43 in ``Viscosity Solutions and Applications''  Ed.'s Dolcetta and Lions, 
SLNM {\bf 1660}, Springer Press, New York, 1997.

 \smallskip

\noindent
\item{[\GAR]}   L. G\aa rding, {\sl  An inequality for hyperbolic polynomials},
 J.  Math.  Mech. {\bf 8}   no. 2 (1959),   957-965.

 \smallskip

 \noindent
\item{[\HAR]}   F. R. Harvey,   {\sl Removable singularities and structure theorems for positive currents}. Partial differential equations (Proc. Sympos. Pure Math., Vol. XXIII, Univ. California, Berkeley, Calif., 1971), pp. 129-133. Amer. Math. Soc., Providence, R.I., 1973.

\smallskip

 \noindent
\item{[\DDD]}   F. R. Harvey and H. B. Lawson, Jr., {\sl  Dirichlet duality and the non-linear Dirichlet problem},    Comm. on Pure and Applied Math. {\bf 62} (2009), 396-443. ArXiv:math.0710.3991.

\smallskip

 \noindent
\item{[\PUP]}  \ \----------,  {\sl  Plurisubharmonicity in a general geometric context},  Geometry and Analysis {\bf 1} (2010), 363-401. ArXiv:0804.1316.

\smallskip

 \noindent
\item{[\DDR]}  \ \----------, {\sl Dirichlet Duality and the Nonlinear Dirichlet Problem on Riemannian Manifolds},  J. Diff. Geom. {\bf 88} (2011), 395-482.   ArXiv:0912.5220.

\smallskip

 \noindent
\item {[\HYP]} \ \----------, {\sl  Hyperbolic polynomials and the Dirichlet problem},   ArXiv:0912.5220.

\smallskip

 \noindent
\item{[\HLGAR]}  \ \----------, {\sl  G\aa rding's theory of hyperbolic polynomials},
   {Communications in Pure and Applied Mathematics}   {\bf 66} no. 7 (2013), 1102-1128.

 \smallskip

 \noindent
\item{[\REST]}  \ \----------, 
{\sl  The restriction theorem for fully nonlinear subequations},   
 {Ann. Inst. Fourier}  (to appear). ArXiv:1101.4850.

\smallskip

 \noindent
\item{[\AC]}  \ \----------,  {\sl  Potential theory on almost complex manifolds},  {Ann. Inst. Fourier}  (to appear).
ArXiv: 1107.2584.
\smallskip

 \noindent
\item{[\pPSH]}  \ \----------,  {\sl  p-convexity, p-plurisubharmonicity  and the Levi problem },
 Indiana Univ. Math. J.  {\bf 62} no. 1 (2014), 149-170.  ArXiv:1111.3895.

\smallskip

% \noindent
%\item  {[HL$_9$]} \ \----------, {\sl  Geometric plurisubharmonicity and convexity - an introduction}.   
%ArXiv:1111.3875

%\smallskip

  \noindent
\item{[\SURVEY]}  \ \----------, {\sl  Existence, uniqueness and removable singularities
for nonlinear partial differential equations in geometry},    pp. 102-156 in ``Surveys in Differential Geometry 2013'', vol. 18,  
H.-D. Cao and S.-T. Yau eds., International Press, Somerville, MA, 2013.
ArXiv:1303.1117.
\smallskip

 \noindent
\item{[\BELL]} \ \----------,  {\sl The equivalence of  viscosity and distributional
subsolutions for convex subequations -- the strong Bellman principle},  
Bulletin Brazilian Math. Soc. {\bf 44} no. 4 (2013),  621-652.  ArXiv:1301.4914.

\smallskip

 \noindent
\item{[\ASPECTS]} \ \----------,  {\sl  Tangents to subsolutions -- existence and uniqueness},  Stony Brook Preprint.

\smallskip

 \noindent
\item{[\HP]}   F. R. Harvey and J.  Polking,  {\sl Removable singularities of solutions of linear partial differential equations},  Acta Math. {\bf 125} (1970), 39 - 56.

\smallskip

 \noindent
\item{[\HPP]}   \ \----------,   {\sl
Extending analytic objects}, Comm. Pure Appl. Math. {\bf 28} (1975), 701-727. 

\smallskip

 \noindent
\item{[\LAB]}  D.Labutin,
{\sl Isolated singularities for fully nonlinear elliptic equations},  J. Differential Equations  {\bf 177} (2001), No. 1, 49-76.

 \smallskip
 
 \noindent
\item{[\LLAB]}  \ \----------, {\sl Singularities of viscosity solutions of fully nonlinear elliptic equations}, 
Viscosity Solutions of Differential Equations and Related Topics, Ishii ed., RIMS K\^oky\^uroku
No. 1287, Kyoto University, Kyoto (2002), 45-57

\smallskip

 \noindent
\item{[\LLLAB]}   \ \----------, {\sl Potential estimates for a class of fully nonlinear elliptic equations}, 
Duke Math. J. {\bf 111} No. 1 (2002), 1-49.

\smallskip

 \noindent
\item{[\LAN]}   N. S.  Landkof,   {Foundations of Modern Potential Theory},  Springer-Verlag, New York, 1972.

\smallskip

 \noindent
\item{[\LEL]}  
P. Lelong, 
Fonctions plurisousharmoniques et formes diffŽrentielles positives,  Gordon and  Breach, Paris-London-New York (Distributed by Dunod Žditeur, Paris) 1968.

\smallskip

 \noindent
\item{[\POGA]}  
A. V. Pogorelov,  {\sl On the regularity of generalized solutions of the equation} 
det $(\partial^2 u/\partial x_i \partial x_j) = \phi(x_1, ... , x_n)>0$, Dokl. Akad. Nauk SSSR 200, 1971, pp. 534Ð537. 

\smallskip

 \noindent
\item{[\POGB]}  
  \ \----------,     {\sl The Dirichlet problem for the n-dimensional analogue of the Monge-Ampre equation}, Dokl. Akad. Nauk SSSR 201, 1971, pp. 790Ð793.

\smallskip

 \noindent
\item{[\SHF]}  
B. Shiffman,  {\sl Extension of positive line bundles and meromorphic maps.},  Invent. Math. {\bf 15} (1972), no. 4, 332-347.

\smallskip

 \noindent
\item{[\TWA]}  N. Trudinger and X.-J. Wang, {\sl Hessian Measures I}, Topol. Meth. in Nonlinear Anal. 
{\bf 10} (1997), 225-239.

\smallskip

 \noindent
\item{[\TWB]}    \ \----------, {\sl Hessian Measures II}, Annals of Math. {\bf 150} no. 2 (1999), 579-604.

\smallskip

 \noindent
\item{[\TWC]}    \ \----------,  {\sl Hessian Measures III}, Journal of Functional Analysis {\bf 193} (2002), 1Ð23.

\smallskip

 \noindent
\item{[\VIT]}   A. Vitolo,  {\sl Removable singularities for degenerate elliptic equations without
conditions on the growth of the solution}, Univ. of Salerno Preprint, Feb. 2014.

\end{document}